\documentclass[psamsfonts]{amsart}

\usepackage{amssymb, amscd, graphicx} 
\usepackage[all]{xy}

\newcommand{\g}{\mathfrak{g}}
\newcommand{\C}{\mathbb{C}}
\newcommand{\tr}[1]{\mathrm{tr}\left(#1\right)}
\newcommand{\ab}{\mathbf{A}}
\newcommand{\xb}{\mathbf{X}}
\newcommand{\yb}{\mathbf{Y}}
\newcommand{\zb}{\mathbf{Z}}
\newcommand{\bb}{\mathbf{B}}

\newcommand{\ub}{\mathbf{U}}

\newcommand{\wb}{\mathbf{W}}

\newcommand{\aq}{/\!\!/}
\newcommand{\X}{\mathfrak{X}}
\newcommand{\R}{\mathfrak{R}}
\newcommand{\hm}{\mathrm{Hom}}
\newcommand{\SL}{\mathrm{SL}(3,\C)}
\newcommand{\ot}{\mathrm{Out}}
\newcommand{\F}{\mathtt{F}}
\newcommand{\G}{\mathfrak{G}}
\newcommand{\xt}{\mathtt{x}}
\newcommand{\yt}{\mathtt{y}}
\newcommand{\wt}{\mathtt{w}}
\newcommand{\bt}{\mathtt{b}}
\newcommand{\rt}{\mathtt{r}}
\newcommand{\Z}{Z_\rt}

\newcommand{\id}{\mathbb{I}}

\newcommand{\ti}[1]{t_{(#1)}}

\newcommand{\Ad}{\mathrm{Ad}}
\newcommand{\gAd}{\mathfrak{g}_{\Ad_\rho}}
\newcommand{\Cn}[1]{C^{#1}(\F_r;\mathfrak{g}_{\mathrm{Ad}_\rho})}
\newcommand{\Hp}{H^1_{\mathrm{par}}(\F_r,\gAd)}
\newcommand{\Zp}{Z^1_{\mathrm{par}}(\F_r,\gAd)}

\newtheorem{theorem}{Theorem}
\newtheorem{lemma}[theorem]{Lemma}
\newtheorem{corollary}[theorem]{Corollary}
\newtheorem{prop}[theorem]{Proposition}

\newtheorem{commnt}[theorem]{Comment}
\newtheorem{definition}[theorem]{Definition}

\title[Poisson Geometry of $\SL$-Character Varieties]{Poisson Geometry of $\SL$-Character Varieties Relative to a Surface with Boundary}
\author[S. Lawton]{Sean Lawton}
\address{Mathematics Department,Instituto Superior T\'ecnico, Lisbon, Portugal}
\email{slawton@math.ist.utl.pt}
\urladdr{http://www.math.ist.utl.pt/$\sim$slawton}
\date{\today}
\subjclass[2000]{Primary 58D29; Secondary 14D20}

\begin{document}

\maketitle

\begin{abstract}
The $\SL$-representation variety $\R$ of a free group $\F_r$ arises naturally by considering surface group representations
for a surface with boundary.  There is a $\SL$-action on the coordinate ring of $\R$ by conjugation.  The geometric points of the subring of invariants of this action is an affine
variety $\X$.  The points of $\X$ parametrize isomorphism classes of completely reducible representations.  We show the coordinate ring $\C[\X]$ is a complex Poisson algebra with respect to 
a presentation of $\F_r$ imposed by the surface.  Lastly, we work out the bracket on all generators when the surface is a three-holed sphere or a one-holed torus.  
\end{abstract}


\section{Introduction}

Surface group representations arise naturally in many contexts:  moduli of geometric structures on surfaces \cite{G1,G2}, knot and link invariants \cite{Si}, spin networks \cite{LP}, and moduli of 
flat principal bundles over a surface \cite{GM, GHJW}.  For surfaces with non-empty boundary and representations of their fundamental groups into subvarieties of $\mathrm{GL}_m$, the space of representations is an affine variety.  The moduli of conjugacy classes of these representations (the categorical quotient) is likewise an affine variety.  The algebraic structure of these varieties cannot alone distinguish the moduli for differing surfaces with the same Euler characteristic.  However, there is additional structure making the moduli a Poisson variety.  This additional structure distinguishes the moduli of surface group representations for surfaces with like Euler characteristic which are not homeomorphic.  We demonstrate this for surfaces with $\chi=-1$.  In general, if two non-homeomorphic surfaces have the same Euler characteristic then they cannot have the same number of boundary components, and we will show that the boundary distinguishes the Poisson structure.

There has been much work in this area in recent years.  A brief, albeit incomplete, history is in order.  In the early $1980$'s, Goldman and Wolpert showed that there is a natural symplectic 
geometry on the space of surface group representations for a surface without boundary \cite{G4, W}.  In Goldman's paper he proves that the symplectic form is closed using the infinite dimensional 
space of all connections on a principal bundle over the surface.  This prompted Karshon in 1992 \cite{Ka} to introduce an algebraic proof that the 2-form is closed.  Then in $1997$ the first proof 
that Goldman's $2$-form extended to surfaces with boundary was published by Guruprasad, Huebschmann, Jeffrey, and Weinstein \cite{GHJW} which built on work that each had done previously.  
This $2$-form was also derived by Kim \cite{Ki} in 1999 in the context of convex projective structures.  However, the moduli of surface group representations is not symplectic when the surface has 
non-empty boundary.  It is foliated by symplectic subspaces and so is a Poisson variety.  In 1986, following his previous work for closed surfaces, Goldman \cite{G5} described the Poisson bracket on 
the moduli in terms of intersections of cycles in homology.  This work was later generalized by Chas and Sullivan in \cite{CS}.  Although it is reasonable to suppose Goldman's formula remains valid 
for surfaces with boundary, no proof exists in the literature.  Moreover, the particular form of the Poisson structures is largely unexplored, as is the relationship between differing structures on 
the same variety.  The only work to date exploring these structures directly is \cite{G3, G6} in the case of $\mathrm{SL}(2,\C)$ representations.  

This paper is organized as follows.  In Section \ref{repchar} we introduce representation and character varieties (the moduli of representations) and review some of their properties.  In Section 
\ref{sympfol} we show the moduli is symplectically foliated using the important work of \cite{GHJW}.  In Section \ref{poissform} we recall Goldman's bracket formula and 
prove that his formulation extends to surfaces with boundary.  Before computing the first known examples of $\SL$ Poisson brackets, we review in Section \ref{lawtonresults} the structure of the 
character variety for Euler characteristic -1 surfaces, the only non-trivial case known (see \cite{L}).  Then, in Sections \ref{threehole} and \ref{onehole} we come to the main theorems of this paper.  
We explicitly derive the Poisson structures for the three-holed sphere and the one-holed torus, and discuss symmetries in their formulations coming from the outer automorphisms of the surface group.

In particular, we show that with respect to certain symmetry operators coming from the forementioned outer automorphisms the Poisson bi-vector field for the three-holed sphere is given by:
$$\mathfrak{a}_{4,-4}\frac{\partial}{\partial \ti{4}}\land \frac{\partial}{\partial \ti{-4}}+(1-\mathfrak{i})\left(\mathfrak{a}_{4,5}\frac{\partial}{\partial \ti{4}}\land \frac{\partial}{\partial 
\ti{5}}\right)$$
whereas the bi-vector for the one-holed torus is given by: \begin{align*}&\Sigma_{1}\bigg(\mathfrak{a}_{1,2}\frac{\partial}{\partial \ti{1}}\land \frac{\partial}{\partial \ti{2}}\bigg)
+\Sigma_2\left(\mathfrak{a}_{3,4}\frac{\partial}{\partial \ti{3}}\land \frac{\partial}{\partial \ti{4}}\right)\\
&+\frac{1}{2}\Sigma_1\Sigma_{2}\bigg(\mathfrak{a}_{1,3}\frac{\partial}{\partial \ti{1}}\land \frac{\partial}{\partial \ti{3}}
+\mathfrak{a}_{1,-3}\frac{\partial}{\partial \ti{1}}\land \frac{\partial}{\partial \ti{-3}}\bigg),\end{align*}
despite the fact that both are identical as varieties.  In each expression $\mathfrak{a}_{i,j}$ are polynomials explicitly derived in Sections \ref{threehole} and \ref{onehole}.

\subsection*{Acknowledgments}
The author thanks his advisor Bill Goldman for introducing him to this topic and for his continuous encouragement and support.  This work was completed at Park City Mathematics Institute in $2006$.  
The author thanks PCMI for its financial support and for providing a stimulating environment to work.  Additionally, he thanks Hans Boden, Carlos Florentino, and Elisha Peterson for fruitful 
conversations.  Lastly, we thank the referee whose comments helped make this paper more readable and well organized.

\section{Representation $\&$ Character Varieties of a Free Group}\label{repchar}

\subsection{Definition and Properties}
Let $\F_r$ be a rank $r$ free group generated by $\{\xt_1,...,\xt_r\},$ and let $\G$ be the linear algebraic group $\mathrm{SL}(m,\C)$.  The map 
$$\hm(\F_r,\G)\longrightarrow \G^{\times r}$$ defined by sending $$\rho \mapsto (\rho(\xt_1),\rho(\xt_2),...,\rho(\xt_r))$$ is a
bijection.  Since $\G^{\times r}$ is the $r$-fold product of irreducible algebraic sets, $\G^{\times r} \cong \hm(\F_r,\G)$ is an affine
variety.  Moreover, since the product of smooth varieties over $\C$ is smooth, $\hm(\F_r,\G)$ is non-singular.  As such we denote $\hm(\F_r, \G)$ by $\R$ and refer to it as the $\G$-{\it representation 
variety of $\F_r$}.

Let $\C[\R]$ be the coordinate ring of $\R$.  Our preceding remarks imply $\C[\R]\cong \C[\G]^{\otimes r}$.  For $1\leq k\leq r,$ define a 
{\it generic matrix} from the complex polynomial ring in 
$m^2r$ indeterminates by
$$\xb_k=\left(
\begin{array}{cccc}
x^k_{11} & x^k_{12} & \cdots & x^k_{1m}\\
x^k_{21} & x^k_{22}  &\cdots & x^k_{2m}\\
\vdots &\vdots & \ddots & \vdots\\
x^k_{m1} & x^k_{m2} &\cdots &x^k_{mm}\\
\end{array}\right).$$
Let $\Delta$ be the ideal of relations defining $\R$ in $\C[x^k_{ij}\ |\ 1\leq i,j \leq m^2,\ 1\leq k\leq r]$.  Therefore,
$\Delta=(\det(\xb_k)-1\ |\ 1\leq k\leq r)$, and 
$$\C[\R]=\mathbb{C}[x^k_{ij}\ | \ 1\leq i,j \leq m^2,\ 1\leq k\leq r]/\Delta.$$

Let $(\xb_1,\xb_2,...,\xb_r)$ be an $r$-tuple of generic matrices.  An element $f\in \C[\R]$ is a function defined in
terms of such $r$-tuples (unique up to an element in $\Delta$).  There is a $\G$-action on $\C[\R]$ given by diagonal conjugation.  Precisely, for $g\in \G$ $$g\cdot
f(\xb_1,\xb_2,...,\xb_r)= f(g^{-1}\xb_1 g,...,g^{-1}\xb_r g).$$  We note that this is well defined since the generators of $\Delta$ are invariant under this action.

A linear algebraic group $\G$ is called \emph{geometrically reductive} if for any rational representation $\G\to \mathrm{GL}_m$ and any nonzero invariant $v$ there exists an invariant homogeneous 
polynomial $f$ so $f(v)\not= 0$.  It is a fact, proved by M. Nagata, that if $\G$ is a geometrically reductive 
group acting rationally on an affine variety $\mathrm{Spec}_{max}(A)$, then $A^\G$ is finitely generated.   The ``unitary trick'' implies that $\mathrm{SL}(m,\C)$ is geometrically
reductive, and so the subring of invariants  $\C[\R]^{\G}$ of the conjugation action defined above is a finitely generated $\C$-algebra (see \cite{D,P1,R}).  

The {\it character variety} $$\X=\mathrm{Spec}_{max}(\C[\R]^{\G})$$ is the irreducible algebraic set whose coordinate ring is the
ring of invariants.  Therefore, $\C[\X]$ includes all polynomial maps of the form $\tr{\xb_{i_1}\xb_{i_2}\cdots\xb_{i_k}},$ where $1\leq i_j\leq r$.
However, it is a consequence of results of Procesi \cite{P1} that all elements of the ring of invariants are polynomials in these functions.  For $\G=\mathrm{SL}(m,\C)$, the maximal value of $k$ 
necessary in the expression $\tr{\xb_{i_1}\xb_{i_2}\cdots\xb_{i_k}}$ is known to be between $m(m+1)/2$ and $m^2$; and for $m=1,2,3,4$ it is known to be equal to the lower bound (see \cite{DF}).

For $r>1$, the Krull dimension of $\X$ is $(m^2-1)(r-1)$ since generic elements have zero dimensional isotropy (see \cite{D}, page $98$).  

There is a regular map $\R\stackrel{\pi}{\to}\X$ which factors through $\R/\G$:  let $\mathfrak{m}$ be a maximal ideal
corresponding to a point in $\R$, then the composite isomorphism $\C\to \C[\R]\to \C[\R]/\mathfrak{m}$ implies that the composite map $\C\to
\C[\R]^\G\to \C[\R]^\G/(\mathfrak{m}\cap \C[\R]^\G)$ is an isomorphism as well.  Hence the contraction $\mathfrak{m}\cap \C[\R]^\G$ is maximal, and since for any $g\in \G$, 
$\left(g\mathfrak{m}g^{-1}\right)\cap \C[\R]^\G=\mathfrak{m}\cap \C[\R]^\G$, $\pi$ factors through $\R/\G$ (see \cite{E}, page $38$).  
Although $\R/\G$ is not generally an algebraic set, $\X$ is the categorical quotient $\R\aq\G$, and since $\G$ is a (geometrically) reductive algebraic group $\pi$ is surjective 
and maps closed $\G$-orbits to points (see \cite{D}).

\subsection{Completely Reducible Representations}
For every representation $\rho\in \R$, $\C^m$ is a $\F_r$-module induced by $\rho$.  A {\it completely reducible} representation is one that is a direct sum of irreducible 
subrepresentations.  Such representations induce a {\it semi-simple} module structure on $\C^m$, and irreducible representations respectively result in {\it simple} modules.  
For any composition series of the $\F_r$-module associated to $\rho$, \begin{equation}\label{compseries} \C^m=V_0\supset V_1\supset\cdots \supset V_{l}=0,\end{equation} there is a semi-simple $\F_r$-module $W=\bigoplus V_{i}/V_{i+1}$.  With respect to a chosen basis of $W$, there exists a completely reducible representation $\rho^{(s)}$.  Its conjugacy class is independent of any 
basis and moreover the Jordan-H\"{o}lder theorem implies that this class is also independent of the composition series.  

We characterize these representations by their orbits when $m=3$ (the more general case is analogous) with the following proposition.

\begin{prop}
Let $\G\rho$ be the conjugation orbit of $\rho\in \R$.  Then $\G\rho$ is closed if and only if $\rho$ is completely reducible.
\end{prop}
\begin{proof}
If $\rho$ is not completely reducible ($l>1$ in \eqref{compseries}), 
then it is reducible.  It is either reducible by a dimension 1 or dimension 2 subspace.    If it is reducible by a dimension 1 subspace then there exists a basis so for any $\wt\in\F_r$ 
it has the form:  $$\left[\begin{array}{ccc}a(\wt)
& b(\wt)& c(\wt)\\0 &d(\wt) &e(\wt)\\0&f(\wt)&g(\wt)\end{array}\right].$$  In this form, conjugating by $$\left[\begin{array}{ccc}1 & 0& 0\\0 &1/n
&0\\0&0&1/n\end{array}\right],$$ and taking the limit as $n\to\infty$ results in $$\left[\begin{array}{ccc}a(\wt)           
& 0& 0\\0 &d(\wt) &e(\wt)\\0&f(\wt)&g(\wt)\end{array}\right].$$  This limiting representation is $\rho^{(s)}$, if it had two irreducible summands.  Otherwise we may conjugate 
$\rho$ so $f(\wt)$ may be taken to be $0$.  Then conjugating this form of $\rho$ by $$\left[\begin{array}{ccc}1 & 0& 0\\0 &1/n
&0\\0&0&1/n^2\end{array}\right]$$ and taking the limit as $n\to\infty$ results in $$\left[\begin{array}{ccc}a(\wt)           
& 0& 0\\0 &d(\wt) &0\\0&0&g(\wt)\end{array}\right],$$ which is $\rho^{(s)}$ when it has three irreducible summands.  Either way, we have a sequence, $\rho_n\in \G\rho$, beginning 
at $\rho$ and limiting to $\rho^{(s)}\notin \G\rho$.  

Now if $\rho$ is reducible by a dimension 2 subspace then there exists a basis so for any $\wt\in\F_r$ $\rho$ has the form:  $$\left[\begin{array}{ccc}a(\wt)
& b(\wt)& c(\wt)\\d(\wt) &e(\wt)&f(\wt)\\0&0&g(\wt)\end{array}\right].$$  In this form, conjugating by $$\left[\begin{array}{ccc}1 & 0& 0\\0 &1
&0\\0&0&1/n\end{array}\right],$$ and taking the limit as $n\to\infty$ results in $$\left[\begin{array}{ccc}a(\wt)
& b(\wt)& 0\\d(\wt) &e(\wt)&0\\0&0&g(\wt)\end{array}\right].$$  This limiting representation is $\rho^{(s)}$, if it had two irreducible summands.  Otherwise we may conjugate
$\rho$ so $d(\wt)$ may be taken to be $0$.  Then conjugating this form of $\rho$ by $$\left[\begin{array}{ccc}1 & 0& 0\\0 &1/n
&0\\0&0&1/n^2\end{array}\right]$$ and taking the limit as $n\to\infty$ results in $$\left[\begin{array}{ccc}a(\wt)
& 0& 0\\0 &e(\wt) &0\\0&0&g(\wt)\end{array}\right],$$ which is $\rho^{(s)}$ when it has three irreducible summands.  Again we have a sequence, $\rho_n\in \G\rho$, beginning
at $\rho$ and limiting to $\rho^{(s)}\notin \G\rho$.

It follows that if $\G\rho$ is closed then $\rho$ is completely reducible.  

For the converse, we first show that $\pi(\rho)=\pi(\psi)$ if and only if $\rho^{(s)}$ is conjugate to $\psi^{(s)}$.  Indeed, suppose that $\pi(\rho)=\pi(\psi)$.  Then their characteristic 
polynomials are equal:  $\chi_{\rho}=\chi_{\psi}$.  Thus $\chi_{\rho^{(s)}}=\chi_{\psi^{(s)}}$.  However semi-simple representations are determined by their characteristic 
polynomials (this is a non-trivial fact established in \cite{A}), so $\rho^{(s)}$ is conjugate to $\psi^{(s)}$.  On the other hand, if for some basis $\rho^{(s)}=\psi^{(s)}$ then $\rho_n\rightarrow\rho^{(s)}=\psi^{(s)}\leftarrow\psi_n$.  This in turn implies $\overline{\G\rho}\cap\overline{\G\psi}$ is not empty, and so $\pi(\rho)=\pi(\psi)$.
  
Now suppose $\rho$ has a non-closed orbit, and let $\psi$ be an element of $\overline{\G\rho}-\G\rho$.  Then $\rho$ and $\psi$ are not conjugate.  Without loss of generality, we 
can assume that $\G\psi$ is closed since the dimension of each subsequent sub-orbit decreases.  So $\psi=\psi^{(s)}$ and $\pi(\psi)=\pi(\rho)$.  Hence, $\rho^{(s)}=\psi$ and so 
$\rho$ cannot be completely reducible else it would be conjugate to $\psi$, which it is not.  In other words, if $\rho$ is completely reducible, then $\G\rho$ is closed.
\end{proof}

Let $\R^{ss}$ be the subset of $\R$ containing only completely reducible representations.  Then we have just shown that $\R^{ss}$ is the set of representations with closed orbits, $\R^{ss}/\G$ is in 
bijective correspondence (as sets) to $\X$, and the following diagram commutes:

\begin{displaymath}
\begin{CD}
\R@>>>\X\\
@AAA      @AAA\\
\R^{ss} @>>> \R^{ss}/\G.
\end{CD}
\end{displaymath}

For a complete treatment of the above arguments see \cite{A,P2}.

\subsection{Simple Representations}

Let $\R^s\subset \R^{ss}$ be the set of irreducible representations, and let $\R^{reg}$ be the regular points in $\R$; that is the representations that have closed orbits
and have minimal dimensional isotropy.  These points form an open dense subset of $\R$ (see \cite{D}).

\begin{prop}
$\R^{reg}=\R^s$, if $\F_r$ has rank greater than $1$.  
\end{prop}

\begin{proof}
We have already seen that the irreducible representations have closed orbits, since they are completely reducible.  So it remains to show that $\rho$ is irreducible if and only if its isotropy has minimal dimension.  First, however, we address the case of $\F_1$.

In this case, all representations have an invariant subspace since the characteristic polynomial always has a root over $\C$.  So there are no irreducible representations, and 
the semi-simple representations $\R^{ss}$ are exactly the diagonalizable matrices since any matrix over $\C$ is conjugate to one in upper-triangular form.  Moreover, the dimension of the isotropy 
of any such matrix is at least $m-1$ since any matrix commutes with itself.  Also the set of matrices with distinct eigenvalues is dense; and any diagonalizable matrix (the representations with 
closed orbits) has a repeated eigenvalue if and only if its isotropy has dimension strictly greater than $m-1$.  So in this case, $\R^{reg}$ is the set of matrices with distinct eigenvalues, and 
$\R^s=\emptyset$.

Otherwise, $\F_r$ has rank at least $2$.  If $\rho\in \R^{ss}$ has an invariant subspace, it has non-zero dimensional isotropy since it fixes at least one line in $\C^m$.  On 
the other hand, the representations that have at least two distinct generic matrices having no shared eigenspaces have isotropy equal to the center, which is generated by the $m^{\mathrm{th}}$
roots of unity and so is zero-dimensional.  If a representation does not have this property then it must be reducible.  Hence, the minimal dimension of isotropy is zero
which is realized if and only if $\rho \in \R^s\subset \R^{ss}$.  Thus when $\F_r$ has $r>1$, then $\R^{reg}=\R^{s}$.  

\end{proof}

In \cite{A}, it is shown that $\R^s/ \G$ is a non-singular variety.  Moreover, in \cite{G2} it is shown that $\G$ acts properly on $\R^s$ (for $m=3$ only, but the argument 
generalizes), and although the action is not effective, the kernel is the center $\mathbb{Z}_m$.  Thus the induced ``infinitesimal'' action on the tangent space is in fact effective, since the tangent map corresponding to the center is zero.  Thus, if $\rho\in\R^s$ the tangent space to an orbit, $T_{\rho}(\mathcal{O}_{\rho})$, is isomorphic to $\mathfrak{g}$, the Lie 
algebra of $\G$.  Together with properness, this implies that $\R^s\to\R\aq \G$ is a local submersion which in turn implies $T_{\rho}(\R\aq \G)\cong T_{\rho}(\R)/T_{\rho}(\mathcal{O}_{\rho})$ whenever 
$\rho$ is irreducible.  

It is not always the case that the tangent space to the quotient is the quotient of tangent spaces.  Just consider representations from the free group of rank $1$ into $\SL$.  The ring of 
invariants is two dimensional and the ring is generated by $\tr{\xb}$ and $\tr{\xb^{-1}}$.  So the ideal is zero and the ring is free.  Consequently it is smooth and the representation sending 
everything to the identity (having maximal isotropy) is a non-singular point.  This illustrates that there can be smooth points in the quotient that have positive-dimensional isotropy.  At these 
points, $T_{\rho}(\R\aq \G)\not\cong T_{\rho}(\R)/T_{\rho}(\mathcal{O}_{\rho}),$ seen by simply comparing dimensions.

\begin{commnt}
In the sense of Geometric Invariant Theory (GIT), the semistable points of the $\G$ action are the entire space $\R$ since for affine reductive $\G$ varieties all points are semistable in general (there are no points for which all invariants vanish).  So the set of semisimple representations is not equal to the set of semistable points.  The semisimple representations are sometimes referred to as polystable in GIT.  However, the stable points are in fact the simple representations.  We make this comment since the $s$ in the notation $\R^s$ can be taken to mean either stable or simple, but the $ss$ in the notation $\R^{ss}$ means only semisimple (not semistable).
\end{commnt}

\section{Symplectic Foliations on Character Varieties}\label{sympfol}

\subsection{The Boundary Map and Generic Leaves}
Let $\Sigma_{n,g}$ be a compact, connected, smooth, orientable surface of genus $g$ with $n>0$ disjoint open disks removed, and let $*$ be a generic point in $\Sigma_{n,g}$.  Its fundamental group has the following presentation:
$$\pi_1(\Sigma_{n,g},*)=\{\xt_1,\yt_1,...,\xt_g,\yt_g,\bt_1,...,\bt_n\ |\ \xt_1\yt_1\xt_1^{-1}\yt_1^{-1}\cdots \xt_g\yt_g\xt_g^{-1}\yt_g^{-1}\bt_1\cdots \bt_n=1\},$$ which is free of
rank $r=2g+n-1$.  And so its Euler characteristic is $\chi(\Sigma_{n,g})=1-r+0=2-2g-n.$  We will refer to such a surface as {\it n-holed}.  So the following is a two-holed genus two surface.

\begin{figure}[h]
\begin{center}
\includegraphics{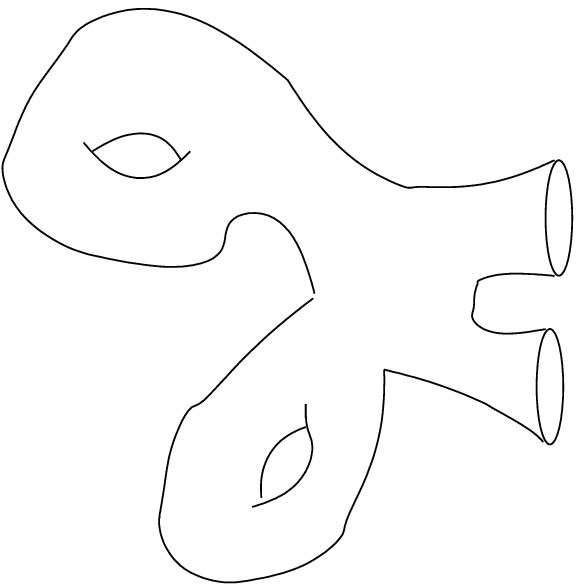}
\caption{$\Sigma_{2,2}$}
\end{center}
\end{figure}

If we assume $r>0$, then $\chi\leq 0$.  In particular, $r=1$ if and only if $\chi=0$, and in this case $\Sigma_{2,0}$ is homeomorphic
to an annulus.  Otherwise, $\chi<0$.  The rank $r$ is 2 if and only if $\chi=-1$, in which case the surface is either
$\Sigma_{3,0}$ or $\Sigma_{1,1}$; that is, the three-holed sphere (or trinion or pair-of-pants), or the one-holed torus,
respectively.

Let $\G=\mathrm{SL}(m,\C)$ and let $\mathbb{I}$ be the $m\times m$ identity matrix.  The coordinate ring of $\G\aq\G$ is freely generated by the $m-1$ coefficients of the characteristic polynomial $$\mathrm{det}(\lambda \mathbb{I}-\xb)=\lambda^m-c_1(\xb)\lambda^{m-1}+\cdots \pm c_{m-1}(\xb)\lambda \mp 1$$ in $\C[\G]$.  Consequently, $\G\aq\G=\C^{m-1}$ which we parametrize by coordinates $(\tau_{(1)},..., \tau_{(m-1)})$.  We note that $c_1(\xb)=\tr{\xb}$ and $c_{m-1}(\xb)=\tr{\xb^{-1}}$ in general.  

We make the following definition in these terms.

\begin{definition}
The $i^{\text{th}}$ boundary map $$\mathfrak{b}_i:\X=\R\aq\G=\hm(\pi_1(\Sigma_{n,g},*),\G)\aq\G\longrightarrow\G\aq\G$$ is the mapping $[\rho]\mapsto [\rho_{|_{\bt_i}}]=(\tau^i_{(1)},..., \tau^i_{(m-1)})$ which sends the extended class of a representation to the class corresponding to the restriction of $\rho$ to the boundary $\bt_i$.
\end{definition}

This is well defined since there exists a representative $\rho$ of $[\rho]$ whose orbit is of maximal dimension (necessarily unique) and for any $g\in \G$ the orbit of $\mathbf{B}_i=\rho(\bt_i)$ and $g\rho(\bt_i)g^{-1}$ are identical.  However, the map $\mathfrak{b}_i$ is invariant and polynomial and any other representative from $[\rho]$ is in the closure of the orbit $\G\rho$.  Hence it is well defined for \emph{any} representative of $[\rho]$.

Subsequently, we make the following further definition.

\begin{definition}
The boundary map is the mapping $\mathfrak{b}=(\mathfrak{b}_1,...,\mathfrak{b}_n):\X=\G^{\times r}\aq\G\longrightarrow (\G\aq\G)^{\times n}$.
\end{definition}

The boundary map $\mathfrak{b}$ depends on the surface, in particular the presentation of its fundamental group.  We refer to it as a {\it peripheral structure}, and the pair
$(\X,\mathfrak{b})$ as the {\it relative character variety}.  

Let $$\mathfrak{F}=\bigcap_{i=1}^{n}\mathfrak{b}_i^{-1}(\tau^i_{(1)},..., \tau^i_{(m-1)})=\{[\rho]\ |\
\mathfrak{b}([\rho])=((\tau^1_{(1)},..., \tau^1_{(m-1)}),...,(\tau^n_{(1)},..., \tau^n_{(m-1)}))\}.$$  Each $\mathfrak{F}$ is cut out of $\X$ by the
equations $c_k(\mathbf{B}_i)=\tau^i_{(k)}$ for $1\leq k\leq m-1$ and $1\leq i\leq n$, and so is an algebraic set.  Moreover, they partition $\X$
since every representation has well-defined boundary values.

The dimension of $\X$ is $(m^2-1)(r-1)=2(m^2-1)(g-1)+(m^2-1)n$, and imposing boundary values provide $(m-1)n$ relations (not necessarily independent) in $\C[\X]$.  Therefore, the dimension of any $\mathfrak{F}$ is greater than or equal to $2(m^2-1)(g-1)+m(m-1)n$.

We now prove the following proposition which we need to show that this partitioning of $\X$ is in fact a foliation of its smooth stratum.

\begin{prop}
If $g\geq 1$ then $\mathfrak{b}$ is surjective.  In all other cases, $\mathfrak{b}$ is dominant; that is $\overline{\mathfrak{b}(\X)}=\C^{n(m-1)}$.
\end{prop}

\begin{proof}
We first suppose $g\geq 1$.  The mapping $\mathrm{GL}(m,\C)\times \mathrm{GL}(m,\C)\to \mathrm{SL}(m,\C)$ given by $(\ab,\bb)\mapsto [\ab,\bb]=\ab\bb\ab^{-1}\bb^{-1}$ is surjective (see \cite{Fo}).  This easily implies that the mapping $\kappa_1:\G\times \G\to \G$ is surjective by simply choosing $m^{\text th}$ roots of the determinants of $\ab$ and $\bb$.  In fact, we have that the mapping $\kappa_g:\G^{\times 2g}\to \G$ given by $(\ab_1,\bb_1,...,\ab_g,\bb_g)\mapsto [\ab_1,\bb_1]\cdots [\ab_g,\bb_g]$ is also surjective simply by letting $\ab_2=\bb_2=\cdots \ab_g=\bb_g=\id$ and picking a value for $\kappa_1$.

Now take any point $p=\left((\tau^1_{(1)},..., \tau^1_{(m-1)}),...,(\tau^n_{(1)},..., \tau^n_{(m-1)})\right)$ in $\C^{n(m-1)}$.  We wish to find $[\rho]\in \X$ so $\mathfrak{b}([\rho])=p$.  Let $r=2g+n-1$, and let $\R=\G^{\times r}$ be the subset of $\R\times \G$ parametrized by $$\left(\xb_1,\yb_1,...,\xb_g,\yb_g,\bb_1,...,\bb_{n-1},(\Pi_{i=1}^g[\xb_i,\yb_i]\Pi_{j=1}^{n-1} \bb_j)^{-1}\right).$$  

For all $1\leq k\leq n$ and every tuple $(\tau^k_{(1)},..., \tau^k_{(m-1)})$ there exists a matrix $\bb_k$ so $c_j(\bb_k)=\tau^k_j$ for $1\leq j\leq m-1$ since there exists a global section to $\G\to \G\aq \G$ given by the rational canonical form of a matrix.  Now let $\zb=\left(\bb_1\cdots \bb_n\right)^{-1}$.  Then there exists $q=(\xb_1,\yb_1,...,\xb_g,\yb_g)$ so $\kappa_g(q)=\zb$.  Then it follows that $\rho:=(q,\bb_1,...,\bb_n)$ is in $\R\subset \R\times \G$.  Letting $[\rho]=\pi(\rho)\in \X$, we have that $\mathfrak{b}([\rho])=p$.

Now if $g=0$ the mapping may not be surjective (this is related to the Deligne-Simpson problem).  However, it is dominant; that is its image is dense.  

Let $\mathcal{O}$ be the subset of $\G$ defined by the property that all eigenvalues are distinct and no proper subset of the eigenvalues multiplies to 1.  Foth \cite{Fo} calls the latter of these two conditions {\it property P} and Kostov \cite{Ko} says such matrices have {\it generic eigenvalues}.  Both conditions define Zariski open sets and so the intersection of these two sets is likewise Zariski open.  Hence, the $m-1$ coefficients of the characteristic polynomial, $\{c_1,...,c_{m-1}\}$, are independent on this open dense set since they are polynomial and independent globally.  Consequently, the image of the $n(m-1)$-tuple $((c_1,...,c_{m-1}),...,(c_1,...,c_{m-1}))$ from $\mathcal{O}^n$ to $\C^{n(m-1)}$ is dense.  Call the image $c(\mathcal{O})$.  Take any point $p\in c(\mathcal{O})$.  We then have $n$ conjugacy classes in $\G$:  $C_1,...,C_n$; where each has representatives with distinct and generic eigenvalues.  It is shown in \cite{Ko} that the set $$\mathcal{W}(C_1,...,C_n):=\{(M_1,...,M_n)\ |\ M_i\in C_i,\ M_1\cdots M_n=\mathbb{I}\}$$ is a smooth connected variety of positive dimension $\sum \mathrm{dim}C_i -(m^2-1)$.  Now take any such tuple of $n$ matrices, then it defines a representation $\rho:\pi_1(\Sigma_{0,n},*)\to \G$ and so $\mathfrak{b}([\rho])=p$.  Hence $c(\mathcal{O})\subset \mathrm{Image}(\mathfrak{b})$.  Since $\overline{c(\mathcal{O})}=\C^{n(m-1)}$, $\mathfrak{b}$ is dominant.

\end{proof}

Consequently, we have the following corollary (see \cite{Ha} page 271).

\begin{corollary}
There exists an open dense subset (Zariski open) $\X^{\mathfrak{b}}\subset \X$ where $\mathfrak{b}$ is a submersion.
\end{corollary}

Let $\mathcal{X}$ be the complement of the singular locus (a proper sub-variety) in $\X$, so $\mathcal{X}$ is a non-singular complex manifold that is dense in $\X$.  The construction of $\X^{\mathfrak{b}}$ provides that $\X^{\mathfrak{b}}\subset \mathcal{X}$.  Thus the implicit function theorem (see \cite{Wa}) on the smooth complex manifold $\X^{\mathfrak{b}}$ gives that $\mathcal{F}:=\mathfrak{F}\cap\X^{\mathfrak{b}}$ is a non-singular complex submanifold of dimension $(m^2-1)(r-1)-(m-1)n=2(m^2-1)(g-1)+m(m-1)n$.  In particular, these submanifolds are of even complex dimension.

Since $\mathfrak{b}$ is a submersion here, the union of these {\it leaves}, $\mathcal{F}$, foliate $\mathfrak{X}^{\mathfrak{b}}$.  We will show that these leaves are complex symplectic submanifolds, which extend to make $\mathcal{X}$ a complex Poisson manifold.  Moreover, the Poisson structure extends continuously over singularities in $\X$.  With respect to this structure, we will refer to the relative character variety $(\X,\mathfrak{b})$ as a {\it Poisson variety}; in particular, an affine variety whose coordinate ring is a Poisson algebra.  We will often drop the fact that we are working over $\C$, but whenever we are discussing a symplectic structure we are referring to a closed $(2,0)$-form.

\subsection{Tangents, Cocycles, and Coboundaries}

Let $\F_r$ be a free group of rank $r$, and let $\mathfrak{g}$ be the Lie algebra of $\G$ identified with its right invariant
vector fields.  For $\rho\in \R$, $\mathfrak{g}$ is a $\F_r$-module, $\mathfrak{g}_{\mathrm{Ad}_\rho}$, given by:
$$\begin{CD}
\F_r  @>\rho>> \G @>\mathrm{Ad}>>\mathrm{Aut}(\mathfrak{g})
\end{CD},$$ where $\mathrm{Ad}(\rho(\wt))(X)=\rho(\wt)X\rho(\wt)^{-1}$ is the adjoint representation.

Define $\Cn{0}=\mathfrak{g}$ and $\Cn{n}=\{\F_r^{\times n}\to \mathfrak{g}\}$, the vector space of functions $\F_r^{\times
n}\to \mathfrak{g}$.  Now define $\delta_n:\Cn{n}\to\Cn{n+1}$ by
\begin{align*}
\delta_n f(\wt_1,\wt_2,...,\wt_{n+1})=&\mathrm{Ad}_\rho(\wt_1)f(\wt_2,...,\wt_{n+1})+\sum_{i=1}^n(-1)^i f(\wt_1,...,\wt_i\wt_{i+1},...,\wt_{n+1})\\
&+(-1)^{n+1}f(\wt_1,...,\wt_n).
\end{align*}

One may verify that $\delta_{n+1}\circ\delta_{n}=0$, and so $(\Cn{*},\delta_*)$ is a cochain complex with coboundary operator
$\delta_*$.

Let $\G^{\F_r^{\times n}}$ be the set of functions $\F_r^{\times n}\to \G$.  Let $f_t$ be a curve in $\G^{\F_r^{\times n}}$, and for
$(\wt_1,...,\wt_n)\in\F_r^{\times n}$ let $\epsilon_{(\wt_1,...,\wt_n)}(f_t)=f_t(\wt_1,...,\wt_n)$ be the evaluation function.  Then for each evaluation, we have a
curve in $\G$, and so we say $f_t$ is {\it smooth} if and only if it is smooth at all evaluations (see \cite{Ka}).  We define the tangent space at a function $f$
to be the vector space of tangents to smooth curves $f_t$ where $f_0=f$.  In other words, $$T_{f}(\G^{\F_r^{\times n}})\cong\mathfrak{g}^{\F_r^{\times
n}}=\Cn{n},$$ given by
$$u_{(\wt_1,...,\wt_n)}(f)=\frac{\mathrm{d}}{\mathrm{d}t}\bigg|_{t=0}\!\!\!\!\epsilon_{(\wt_1,...,\wt_n)}(f_t)=
\frac{\mathrm{d}}{\mathrm{d}t}\bigg|_{t=0}\!\!\!\mathrm{exp}(t\alpha_{(\wt_1,...,\wt_n)})\epsilon_{(\wt_1,...,\wt_n)}f,$$ where $\alpha\in \mathfrak{g}^{\F_r^{\times
n}}$.  Since a function is determined by its evaluations, we consider right invariant vector fields defined along these coordinates which give paths with the requisite
properties.

Let $I$ be a finite subset of $\F_r$, and let $\C[\G^I]$ be the coordinate ring of $\G^I$.  The set of such $I$'s is partially ordered by set inclusion and so
$\C[\G^J]\hookrightarrow \C[\G^I]$ for $J\subset I$.  We then define $\C[\G^{\F_r}]=\underrightarrow{\lim}\C[\G^I]$.  With this said, we note that $\R$ is the subspace of
$\G^{\F_r}$ that is cut out by the functions $\mathtt{ob}_{\xt,\yt}(f)=f(\xt)f(\yt)f(\xt\yt)^{-1}$.  In other words, letting $\id$ be the $m \times m$ 
identity matrix, $$\C[\R]=\C[\G^{\F_r}]/(\mathtt{ob}_{\xt,\yt}-\id\ |\ \xt,\yt\in \F_r),$$ and $\R=\mathrm{Spec}_{max}(\C[\R])$.

The following proposition is a classical fact which we prove here to keep the discussion self-contained.

\begin{prop}
If $\rho \in \R^{reg}=\R^{s}\subset \R^{ss}$, and $[\rho]=\pi(\rho)\in \X$, then
$$T_{[\rho]}(\X)\cong
T_\rho(\R)/T_\rho(\mathcal{O}_\rho)=Z^1(\F_r;\mathfrak{g}_{\mathrm{Ad}_\rho})/B^1(\F_r;\mathfrak{g}_{\mathrm{Ad}_\rho})=H^1(\F_r;\mathfrak{g}_{\mathrm{Ad}_\rho}).$$
\end{prop}

\begin{proof}
In the terms outlined above, $T_f(\R)=\{u\in T_f(\G^{\F_r})\ | \ u(f(\xt)f(\yt)f(\xt\yt)^{-1})=0\}$, and so
\begin{align*}
0&=u(f(\xt)f(\yt)f(\xt\yt)^{-1})\\
&=u(f(\xt))f(\yt)f(\xt\yt)^{-1}+f(\xt)u(f(\yt))f(\xt\yt)^{-1}-f(\xt)f(\yt)f(\xt\yt)^{-1}u(f(\xt\yt))f(\xt\yt)^{-1},
\end{align*}
which implies $u_{\xt\yt}=u_\xt+\mathrm{Ad}_f(\xt)u_\yt.$

However, this is exactly the condition for $\delta_1(f)=0$, and so $$T_f(\R)=\mathrm{Ker}(\delta_1)=Z^1(\F_r;\mathfrak{g}_{\mathrm{Ad}_f}).$$

On the other hand, consider a smooth path contained in the orbit $\mathcal{O}_f\subset \R$:
$$f_t(\xt)=\mathrm{exp}(tu_\xt)f(\xt)=\mathrm{exp}(-tu_0)f(\xt)\mathrm{exp}(tu_0),$$ for
some $u_0\in \mathfrak{g}=\Cn{0}$.  Then
$$u_\xt f(\xt)=\frac{\mathrm{d}}{\mathrm{d}t}\bigg|_{t=0}\!\!\!\!f_t(\xt)=-u_0f(\xt)+f(\xt)u_0,$$ which implies
$u_\xt=\mathrm{Ad}_f(\xt)u_0-u_0$.  However this is exactly the condition $\delta_0(u_0)=u$, so
$$T_f(\mathcal{O}_f)=\mathrm{Image}(\delta_0)=B^1(\F_r;\mathfrak{g}_{\mathrm{Ad}_f}).$$

Therefore, if $\rho \in \R^{reg}=\R^{s}\subset \R^{ss}$, and $[\rho]=\pi(\rho)\in \X$, then
$$T_{[\rho]}(\X)\cong
T_\rho(\R)/T_\rho(\mathcal{O}_\rho)=Z^1(\F_r;\mathfrak{g}_{\mathrm{Ad}_\rho})/B^1(\F_r;\mathfrak{g}_{\mathrm{Ad}_\rho})=H^1(\F_r;\mathfrak{g}_{\mathrm{Ad}_\rho}),$$ as was to be shown.
\end{proof}

\subsection{Homology and the Fundamental Cycle}

Let $\mathbb{Z}\F_r$ be the integral group ring of $\F_r$ and let $\epsilon:\mathbb{Z}\F_r\to \mathbb{Z}$ be the {\it
augmentation map} defined by $\sum n_{\wt}\wt\mapsto \sum n_{\wt}$.  Define $\mathbb{Z}$-modules
$C_0(\F_r)=\mathbb{Z}$ and $C_n(\F_r)=\mathbb{Z}\F_r^{\times n}$, and let $$\partial_{n+1}:C_{n+1}(\F_r)\longrightarrow
C_{n}(\F_r)$$ be defined by
\begin{align*}
\partial_{n+1}(\wt_1,\wt_2,...,\wt_{n+1})=&\epsilon(\wt_1)(\wt_2,...,\wt_{n+1})+\sum_{i=1}^{n}(-1)^i(\wt_1,...,\wt_i\wt_{i+1},...,\wt_{n+1})\\
&+(-1)^{n+1}(\wt_1,...,\wt_n)\epsilon(\wt_{n+1}).
\end{align*}

One can show $\partial_{n}\circ\partial_{n+1}=0$ and so $(C_*(\F_r),\partial_*)$ is a chain complex with boundary $\partial_*$.

In \cite{F} it is shown that the derivations on $\mathbb{Z}\F_r$, $$\mathrm{Der}(\F_r)=\{ D:\mathbb{Z} \F_r \to
\mathbb{Z} \F_r \ | \  D( \xt \yt ) = D ( \xt ) \epsilon ( \yt ) + \xt D ( \yt ) \},$$
are freely generated by the derivations $\frac{\partial}{\partial\xt_i}(\xt_j)=\delta_{ij}$.  These derivations
and their generators are called the \emph{Fox derivatives}.

Recall that the peripheral structure on $\X$ is given by the presentation
$$\F_r=\pi_1(\Sigma_{n,g},*)=\{\xt_1,\yt_1,...,\xt_g,\yt_g,\bt_1,...,\bt_n\ |\ \rt=1\},$$
where $$\rt=\xt_1\yt_1\xt_1^{-1}\yt_1^{-1}\cdots \xt_g\yt_g\xt_g^{-1}\yt_g^{-1}\bt_1\cdots \bt_n=\Pi[\xt_i,\yt_i]\Pi \bt_j.$$  Then with respect to the Fox derivatives and this presentation \cite{Ki} 
shows
\begin{align*}
\partial_2 \Z =&\sum_{j=1}^n \bt_j\ \mathrm{where}\\
\Z=&\sum_{i=1}^g(\frac{\partial \rt}{\partial \xt_i},\xt_i)+(\frac{\partial \rt}{\partial
\yt_i},\yt_i)+\sum_{j=1}^n(\frac{\partial \rt}{\partial \bt_j},\bt_j).
\end{align*}

Consequently, he refers to $\Z$ as the \emph{fundamental relative cycle}, since $$\Z-\sum_{j=1}^n(\frac{\partial \rt}{\partial \bt_j},\bt_j)$$ is the fundamental
cycle when $n=0$, and consequently $[\Z]$ is a generator of $H_2(\F_r,\{\bt_1,...,\bt_n\};\mathbb{Z})\cong \mathbb{Z}$ (see \cite{GHJW}).

\subsection{Parabolic Cocycles}

Let $\F^i_r\subset \F_r$ be the cyclic subgroup generated by the boundary curve $\bt_i$.

\begin{definition}
The set of {\it parabolic cocycles}, $Z^1_{\mathrm{par}}(\F_r;\mathfrak{g}_{\mathrm{Ad}_\rho})\subset
Z^1(\F_r;\mathfrak{g}_{\mathrm{Ad}_\rho})$, is defined by
$f\in Z^1_{\mathrm{par}}(\F_r;\mathfrak{g}_{\mathrm{Ad}_\rho})$ if and only if $f_i=f|_{\F_r^i}\in
B^1(\F^i_r,\mathfrak{g}_{\mathrm{Ad}_\rho})$ for all $1\leq i\leq n$.
\end{definition}
It is shown in \cite{Ki} that $$B^1(\F_r;\mathfrak{g}_{\mathrm{Ad}_\rho})\subset
Z^1_{\mathrm{par}}(\F_r;\mathfrak{g}_{\mathrm{Ad}_\rho}) \subset Z^1(\F_r;\mathfrak{g}_{\mathrm{Ad}_\rho}).$$

So for $[\rho] \in \R^{reg}\aq \G=\X^{reg}$, we have $$\Hp=Z^1_{\mathrm{par}}(\F_r;\mathfrak{g}_{\mathrm{Ad}_\rho})
/B^1(\F_r;\mathfrak{g}_{\mathrm{Ad}_\rho})\ \subset T_{[\rho]}(\X).$$

In other words, $\Hp$ is the set of tangents that are zero on the boundary;
that is, the tangents to representations with constant boundary value.  

\begin{prop}\label{foliation}
$\mathcal{X}$ is foliated by the fibers of the boundary map.
\end{prop}

\begin{proof}
The distribution $\Hp\subset T_{[\rho]}(\X^{reg})$ is exactly the tangents to curves in $\mathfrak{F}\cap \X^{reg}$.

On the other hand, let $\rho\mapsto \tr{\rho(\wt)}$ for a fixed word $\wt$ be a ``word map.''  The image of
sufficiently many such maps (necessarily finite) determines $[\rho]$.  Then holding the
boundary values fixed gives functions in $\C[\mathfrak{F}]\subset \C[\X]$.  As the boundary values are deformed to a different leaf, the word map is likewise deformed.
Therefore, the word maps generate a family of smooth invariant vector fields that generate the distribution.

Consequently, the Stefan-Sussmann theorem (see page $17$ in \cite{DZ}) implies that $\X^{reg}\subset \mathcal{X}$ is foliated by $\mathfrak{F}\cap\X^{reg}$.  But since $\X^{reg}$ is an open dense set and all vector fields corresponding to the distribution are globally defined and continuous on $\X$ and non-singular on $\mathcal{X}$ they can be extended to a foliation on $\mathcal{X}$, with leaves given by $\mathcal{X}\cap \mathfrak{F}$.  We note that the leaves from the submersive points of $\mathfrak{b}$, $\mathcal{F}$, are dense in $\mathcal{X}\cap \mathfrak{F}$.

\end{proof}

\subsection{$(2,0)$-Form on Symplectic Leaves}

In \cite{GHJW}, it is shown that $\omega$ from the following commutative diagram defines a symplectic form on $\mathfrak{F}\cap\mathfrak{X}^{reg}$:
$$ \xymatrix{
H^1(\Sigma_{n,g},\partial \Sigma_{n,g};\gAd) \times H^1(\Sigma_{n,g};\gAd) \ar[r]^-\cup   & H^2(\Sigma_{n,g},\partial \Sigma_{n,g};\gAd \otimes \gAd) \ar[d]^{\mathrm{tr}_*}\\
    &  H^2(\Sigma_{n,g},\partial \Sigma_{n,g};\C) \ar[d]^{\cap [\Z]}\\
H^1_{\mathrm{par}}(\Sigma_{n,g};\gAd)\times H^1_{\mathrm{par}}(\Sigma_{n,g};\gAd)\ar[uu] \ar[r]^-\omega  & H_0(\Sigma_{n,g};\C)=\C. }$$

We note that in the above diagram $H^1_{\mathrm{par}}(\Sigma_{n,g};\gAd)$ arises from
$$j\!\!:\! H^1(\Sigma_{n,g},\partial \Sigma_{n,g};\gAd) \longrightarrow H^1(\Sigma_{n,g};\gAd)$$ as the right factor by $$\mathrm{im}(j)=H_{\text{par}}^1(\Sigma_{n,g};\gAd)\subset  H^1(\Sigma_{n,g};\gAd)$$ and as the left 
factor by $$H_{\text{par}}^1(\Sigma_{n,g};\gAd) \cong H^1(\Sigma_{n,g},\partial \Sigma_{n,g};\gAd)/\mathrm{ker}(j).$$  And the map is well defined since, with respect to the cup product, $\ker(j)$ is orthogonal to 
$\mathrm{im}(j)$.  Lastly, the (relative) group cohomology is isomorphic to the singular cohomology since we are working with a $K(\pi,1)$ space.

It is clear from its definition that $\omega$ is well-defined, bilinear, skew-symmetric; and that for all vector fields $\xi_1$ and $\xi_2$ taking values in $\Hp$,
$$\omega(\xi_1,\xi_2): \X^{reg}\to \Hp^{\times 2}\to \C,$$ is an element of $\C[\X]$ and thus smooth.  Moreover, since the trace pairing is non-degenerate $\omega$ is as well.

Consequently, it is a $(2,0)$-form.  For it to be symplectic, one must check that it is \emph{closed}.  This is not at all obvious, and there are many differing proofs. See \cite{GHJW} for a brief 
history, and a self-contained algebraic proof.

Subsequently, using this form with Proposition \ref{foliation} establishes 
\begin{corollary}
The smooth leaves of $\X$ are symplectic and hence $\mathcal{X}$ is a Poisson manifold, and $\X$ is a Poisson variety (as defined at the beginning of this section).  
\end{corollary}
 
Moreover, having the symplectic form in terms of cup and cap products allows one to verify that the Poisson bracket formula derived in \cite{G5} for closed surfaces directly generalizes to 
surfaces with boundary.  In \cite{Ki} this is not possible, since although he defines the $2$-form and verifies all requisite properties to establish that it is symplectic, he does not show it arises 
from cup and cap products.  This is one of the many hallmarks of \cite{GHJW}.  

\begin{commnt}
It is worth noting that the construction of $\omega$ and the proof of its properties in \cite{GHJW} are carried out with respect to the ground field $\mathbb{R}$.  However, both its construction and the proof of it properties (in particular closure) are done algebraically with respect to group cohomology.  Consequently, the ground field may be replaced by $\C$ and, in particular, the proof that the form is closed remains unchanged.
\end{commnt}

\begin{commnt}
The reader will notice that the fact $\G=\mathrm{SL}(m,\C)$ was not used to construct the 2-form; consequently much of the above discussion generalizes to character varieties defined 
over any reductive linear algebraic group.  However, in this generality we lose many of the explicit facts concerning singularities, dimensions, generators, etc. that depend on the choice of $\G$.
\end{commnt}

\section{Goldman's Bracket Formula}\label{poissform}

In this section we prove that Goldman's bracket formula extends to surfaces with boundary.  The important ideas that follow are those of Goldman (see \cite{G5}), we only apply the results of 
\cite{GHJW} to his proof.

Let $\Sigma$ be an compact, connected, oriented surface, and let $\alpha, \beta$ be immersed curves in $\Sigma$ representing elements in $\pi_1(\Sigma)$ (in generic position).  Denote the set of 
(transverse) double point intersections of $\alpha$ and $\beta$ by $\alpha\cap\beta$, and the oriented intersection number of $\alpha$ and $\beta$ at
$p\in \alpha\cap\beta$ by $\epsilon(p,\alpha, \beta)$.  Additionally, let $\alpha_p\in \pi_1(\Sigma,p)$ represent the curve $\alpha$ based at $p$.   

Then \cite{G5} proves the following theorem for {\it closed} surfaces.

\begin{theorem}[Goldman]Let $\C[\X]$ be the coordinate ring of the $\mathrm{SL}(m,\C)$ character variety of a closed surface group.  Then
\begin{align*}
\{\mathrm{tr}(\rho(\alpha)),\mathrm{tr}(\rho(\beta))\}&=\sum_{p\in
\alpha\cap\beta}\epsilon(p,\alpha,\beta)\big((\mathrm{tr}(\rho(\alpha_p\beta_p))-(1/m)\mathrm{tr}
(\rho(\alpha))\mathrm{tr}(\rho(\beta))\big)
\end{align*}
defines a Lie bracket on $\mathbb{C}[\mathfrak{X}]$ that is a derivation; in other words a Poisson bracket.
\end{theorem}

The construction is in terms of the symplectic form for the case of closed surfaces which is likewise in terms of cohomology.  Applying the $(2,0)$-form $\omega$ discussed above to the case of surfaces with boundary, we will see this same formula for the bracket remains valid on the leaves $\mathfrak{F}$.

Heuristically, this is not surprising.  The reason is the bracket is summed over intersections of cycles.  The boundary components are disjoint and so bracket to zero.  Such functions that ``Poisson 
commute'' or bracket to zero with all other functions are called {\it Casimirs}.  In parabolic cohomology, these are eliminated and so we remove the Casimirs.  A bracket with no Casimirs is 
symplectic, so if these are the only Casimirs, the formula is valid.

\subsection{Proof of the Bracket Formula}

It is the purpose of this section to prove

\begin{theorem}\label{bracketboundary}
 Let $\C[\X]$ be the coordinate ring of the $\mathrm{SL}(m,\C)$ character variety of a surface group having non-empty boundary.  Then
\begin{align}\label{bracket}
\{\mathrm{tr}(\rho(\alpha)),\mathrm{tr}(\rho(\beta))\}&=\sum_{p\in
\alpha\cap\beta}\epsilon(p,\alpha,\beta)\big((\mathrm{tr}(\rho(\alpha_p\beta_p))-(1/m)\mathrm{tr}
(\rho(\alpha))\mathrm{tr}(\rho(\beta))\big)
\end{align}
defines a Lie bracket on $\mathbb{C}[\mathfrak{X}]$ that is a derivation; in other words a Poisson bracket.
\end{theorem}

The Poisson bracket of two functions on a symplectic leaf of a Poisson manifold is given by $\{f_\alpha, f_\beta\}=\omega(Hf_\alpha,Hf_\beta)$, where $Hf_\alpha$ is the Hamiltonian vector field 
associated to $f_\alpha$.  We will first describe the Hamiltonians in terms of cycles and then use the fact that the 2-form is in terms of cup and cap products to relate these to intersection pairings.  
This will prove the formula to be valid.

We begin by setting out the terms we will need for the discussion and then stating a lemma that will reduce the proof of Theorem \ref{bracketboundary} to a computation.

Let $f:\G\to \C$ be an invariant function in $\C[\G]$.  Then $df$ is a $1$-form on $\G$.  Let $\mathbf{A}\in \G$, then $df(\mathbf{A})\in T_\mathbf{A}^*(\G)$ and extends to a right invariant 
$1$-form $F^*(\mathbf{A})\in \g^*$. Here $F^*:\G\to \g^*$ is defined (for $X\in\mathfrak{g}$) by $$F^*(\mathbf{A}):X\mapsto \frac{d}{dt}\Big|_{t=0}f(\mathrm{exp}(tX)\mathbf{A}).$$

The trace is a non-degenerate symmetric bilinear form $\g\times\g\to \C$.  Thus as $\G$-modules under the adjoint representation, $\g\cong \g^*$ by $\mathrm{tr}_{*} :X\mapsto
\tr{X\underline{\ \ \ }}$.
Let $F=\mathrm{tr}^{-1}_{*} F^*:\G\to\g^*\to\g$.  Then for $f(\mathbf{A})=\tr{\mathbf{A}}$ we have
$$\tr{XF(\mathbf{A})}=\mathrm{tr}_*F(\mathbf{A})(X)=F^*(\mathbf{A})(X)=\frac{d}{dt}\Big|_{t=0}f(\mathrm{exp}(tX)\mathbf{A})$$ 
$$=\tr{X\mathbf{A}}=\tr{X\mathbf{A}}-\frac{1}{m}\tr{\mathbf{A}}\tr{X}=\tr{X\left(\mathbf{A}-\frac{1}{m}\tr{\mathbf{A}}\id\right)};$$ and so
because the trace is non-degenerate $$F(\mathbf{A})=\mathbf{A}-\frac{1}{m}\tr{\mathbf{A}}\id.$$

Now consider $\alpha\in \pi_1(\Sigma,*)$, where we choose a representative that is an immersed curve in $\Sigma$ that intersects itself only at transverse double points.  Let $\rho:\pi_1(\Sigma,p)\to\G$
represent $[\rho]\in \X$ and $\alpha_p\in\pi_1(\Sigma,p)$ correspond to $\alpha$ where $p$ is a simple point (not a point of self intersection).

The Hamiltonian vector field of the \emph{Goldman function} $f_\alpha(\rho)=\tr{\rho(\alpha)}$, a map $\X\to \C$, is given implicitly by $df_\alpha=\omega(Hf_\alpha,\underline{\ \ \ })$. In other 
words, $Hf_\alpha=\omega_*^{-1}\left(df_\alpha\right),$ where $\omega_*$ is the isomorphism $$\Hp\to \Hp^*$$ given by $v\mapsto\omega(\underline{\ \ },v)$.  The chain $\alpha_p\otimes 
F(\rho(\alpha_p))$ is a cycle since $\Ad_\rho(\alpha_p)F(\rho(\alpha_p))=F(\rho(\alpha_p))$, and consequently represents a relative cycle in 
$H_1(\F_r,\{\bt_1,...,\bt_n\};\gAd)$ since non-trivial cycles are relative cycles.

The following lemma describes how the Hamiltonian relates to this relative cycle.

\begin{lemma}\label{hamilcycle}
The relative cycle $\alpha_p\otimes F(\rho(\alpha_p))$ is equal to the cap product of the Hamiltonian $Hf_\alpha$ with the fundamental relative cycle $\Z$.  In other words,  $\alpha_p\otimes 
F(\rho(\alpha_p))=Hf_\alpha\cap[Z_r]$.
\end{lemma}

Assuming Lemma \ref{hamilcycle} we now prove Theorem \ref{bracketboundary}.

\begin{proof}[Proof of Theorem \ref{bracketboundary}]
To derive the bracket formula, we will need to use the intersection pairing of cycles.  Since we are working with a $K(\pi,1)$ we may pass from group cohomology 
to singular cohomology, work with intersection theory, and then return to group cohomology.

We quickly remind the reader of facts concerning intersection pairings which may be found on page 367 in \cite{Br}.  Let $D$ be the inverse of the Poincare duality map $$D:H_1(\Sigma;\gAd)\to 
H^1(\Sigma,\partial \Sigma;\gAd);$$ that is, $D(u)\cap [Z_r]=u$.  The intersection pairing $$\bullet: H_1(\Sigma,\partial \Sigma;\gAd)\times H_1(\Sigma,\partial \Sigma;\gAd)\to  H_2(\Sigma,\partial \Sigma;\gAd)$$ is defined by 
$$u\bullet v=(D(v)\cup D(u))\cap [Z_r].$$  It has the geometric consequence of being equal to the sum taken over transverse intersections of 1-cycles with coefficients.  Precisely, 
$$(\alpha\otimes F) \bullet (\beta\otimes F) =\sum_{p\in \alpha \cap \beta} \epsilon(p,\alpha,\beta)F(\alpha_p)\otimes F(\beta_p).$$

We now compute the Poisson bracket of two invariant functions using Lemma \ref{hamilcycle}:
\begin{align*}
\{f_\alpha, f_\beta\}=&\omega(Hf_\alpha,Hf_\beta)\\
                     =&\mathrm{tr}_*(Hf_\alpha\cup Hf_\beta)\cap [Z_r]\\
                     =&\mathrm{tr}_*\left(\left(D(Hf_\alpha\cap[Z_r])\right)\cup \left(D(Hf_\beta\cap[Z_r])\right)\cap [Z_r]\right)\\
                     =&\mathrm{tr}_*\left(\left(Hf_\alpha \cap [Z_r]\right)\bullet \left(Hf_\beta \cap [Z_r]\right)\right)\\
                     =&\mathrm{tr}_*\left(\alpha\otimes F(\rho(\alpha_p))\bullet \beta\otimes F(\rho(\beta_p))\right)\\
                     =&\sum_{p\in \alpha \cap \beta}\epsilon(p,\alpha,\beta)\tr{F(\rho(\alpha_p))F(\rho(\beta_p))}\\
                     =&\sum_{p\in \alpha \cap \beta}\epsilon(p,\alpha,\beta) \tr{\left(\rho(\alpha_p)-(1/m)\tr{\rho(\alpha_p)}\id\right)\left(\rho(\beta_p)-(1/m)\tr{\rho(\beta_p)}\id\right)}\\
                     =&\sum_{p\in \alpha \cap \beta}\epsilon(p,\alpha,\beta) \left(\tr{\rho(\alpha_p\beta_p)}-(1/m)\tr{\rho(\alpha_p)}\tr{\rho(\beta_p)}\right).
\end{align*}

This concludes the proof of the bracket formula.
\end{proof}

It now remains to prove Lemma \ref{hamilcycle}.

\begin{proof}[Proof of Lemma \ref{hamilcycle}]

Let $[u]\in \Hp$, so $u\in\Zp$; and let $\rho_t$ be a path in $\R$ so $u$ is a tangent.  Then $$(df_\alpha)_\rho:[u]\mapsto 
\frac{d}{dt}\Big|_{t=0}f_\alpha(\rho_t(\alpha_p))=F^*(\rho(\alpha_p))(u(\alpha_p))=\tr{F(\rho(\alpha_p))u(\alpha_p)}.$$

We will now define a number of mappings which we will show are related by a commutative diagram; the lemma will follow.

The trace gives a duality pairing $\mathrm{d\!p}:\g\to\g^*$ given by $u\mapsto \tr{u\underline{\ \ }}$ which pushes forward to $$\mathrm{d\!p}_*:H^1(\F_r;\gAd)\to 
H^1(\F_r;\gAd^*).$$   Consequently, $(\mathrm{d\!p}_*)^*:H^1_{\mathrm{par}}(\F_r;\gAd^*)^*\to H^1_{\mathrm{par}}(\F_r;\gAd)^*$ defined by $v\mapsto v(\mathrm{d\!p}_*(\underline{\ \ }))$ is an 
isomorphism.

Using the Kroneker product define $$\mathrm{k\!p}:H_1(\F_r,\{\bt_1,...,\bt_n\};\gAd)\to 
H^1_{\mathrm{par}}(\F_r;\gAd^*)^*$$ by $\alpha\otimes F(\rho(\alpha_p))\mapsto \mathrm{tr}_*\big<\mathrm{d\!p}_*^{-1}(\underline{\ \ }),\alpha\otimes F(\rho(\alpha_p))\big>.$  We claim 
$\mathrm{k\!p}$ is well-defined.  Let $u$ be a parabolic cycle, $\beta$ a 2-chain, and $\gamma$ a 1-chain on the boundary.  Then $$u(\alpha\otimes F+ \partial \beta +\gamma)=u(\alpha\otimes F)+u(\partial 
\beta)+u(\gamma)=u(\alpha\otimes F)+\delta u(\beta) + 0= u(\alpha \otimes F)$$ because $\delta u=0$ since $u$ is a cycle, and $u(\gamma)=0$ since $u$ is parabolic.  This verifies that 
$\mathrm{k\!p}$ is well-defined.

Lastly, in terms of Poincare duality let 
$\mathrm{P\!D}:H^1_{\mathrm{par}}(\F_r;\gAd) \to H_{\mathrm{par}}^1(\F_r;\gAd^*)^*$ be the isomorphism $$u\mapsto\mathrm{tr}_*\left(\mathrm{d\!p}_*^{-1}(\underline{\ \ })\cup u\right)\cap[Z_r].$$

In these terms, we have
\begin{lemma}\label{diagram} The following diagram commutes:
$$ \xymatrix{
H^1_{\mathrm{par}}(\F_r;\gAd) \ar[r]^-{\cap[Z_r]} \ar[d]^{\omega_*} \ar[dr]^{\mathrm{P\!D}} & H_1(\F_r,\{\bt_1,...,\bt_n\};\gAd) \ar[d]^{\mathrm{k\!p}}\\
H_{\mathrm{par}}^1(\F_r;\gAd)^* &  \ar[l]_-{(\mathrm{d\!p}_*)^*} H_{\mathrm{par}}^1(\F_r;\gAd^*)^*} $$
\end{lemma}

We defer proving Lemma \ref{diagram} until after the proof of Lemma \ref{hamilcycle}.

To complete the proof of Lemma \ref{hamilcycle} we begin with a relative cycle $\alpha \otimes F$ starting in the upper right corner of the diagram: 
$$(\mathrm{d\!p}_*)^*(\mathrm{k\!p}(\alpha\otimes F))=(\mathrm{d\!p}_*)^*\mathrm{tr}_*\big<\mathrm{d\!p}^{-1}_*(\underline{\ \ }), \alpha\otimes F\big>=
\mathrm{tr}_*\big<\mathrm{d\!p}^{-1}_*\mathrm{d\!p}_*(\underline{\ \ }),\alpha\otimes F \big>=\mathrm{tr}_*\big<\ \underline{\ \ }\ ,\alpha\otimes F 
\big>=df_\alpha.$$  However, since $(\mathrm{d\!p}_*)^*$ is an isomorphism and $\mathrm{k\!p}$ is injective on the 
image of $\Hp$ under $\cap[\Z]$, we have
$$\alpha_p\otimes F(\rho(\alpha_p))=\mathrm{k\!p}^{-1}\left((\mathrm{d\!p}_*)^*\right)^{-1}(df_\alpha)=(\omega_*)^{-1}(df_\alpha)\cap[Z_r]
=Hf_\alpha\cap[Z_r],$$ by commutativity of the above diagram.  Note that $\mathrm{k\!p}$ is injective on the image of $\Hp$ under $\cap[\Z]$ because $\mathrm{P\!D}$ is an isomorphism and by 
commutativity.  This completes the proof.

\end{proof}

It remains to prove Lemma \ref{diagram}.

\begin{proof}[Proof of Lemma \ref{diagram}]
First we show the lower triangle commutes: $$(\mathrm{d\!p}_*)^*\mathrm{P\!D}(v)=(\mathrm{d\!p}_*)^*\left(\mathrm{tr}_*\left(\mathrm{d\!p}_*^{-1}(\underline{\ \ })\cup 
v\right)\cap[Z_r]\right)=\mathrm{tr}_*\left(\mathrm{d\!p}_*^{-1}\mathrm{d\!p}_*(\underline{\ \ })\cup v\right)\cap [Z_r]=\omega(\underline{\ \ },v)=\omega_*(v).$$

Second, we show the upper triangle commutes (see \cite{B} page 113): $$\mathrm{k\!p}(v\cap[Z_r])=\mathrm{tr}_*\big<\mathrm{d\!p}_*^{-1}(\underline{\ \ }), v\cap 
[Z_r]\big>=\big<\mathrm{tr}_*\left(\mathrm{d\!p}_*^{-1}(\underline{\ \ })\cup 
v\right),[Z_r]\big>=\mathrm{tr}_*\left(\mathrm{d\!p}_*^{-1}(\underline{\ \ })\cup v\right) \cap [Z_r]=\mathrm{P\!D}(v).$$
\end{proof}

\begin{commnt}
The above argument can be generalized to any reductive linear algebraic group.  Such groups have a symmetric, non-degenerate bilinear form $\mathfrak{B}$ on their Lie algebra that is invariant under the 
adjoint representation.  In the above argument this was the trace form.  To generalize the argument replace $\mathrm{tr}_*$ by $\mathfrak{B}_*$ and replace the explicit form of $F$ by its 
general form defined by $\mathfrak{B}(F(\mathbf{A}),X)=\frac{d}{dt}\Big|_{t=0}f(\mathrm{exp}(tX)\mathbf{A})$.  In special cases, like $\G=\mathrm{SL}(m,\C)$, $F$ can be computed explicitly as we did 
above.  See \cite{G5} for further details.
\end{commnt}

\begin{commnt}
Since the bracket formula depends on a choice of orientation of the surface, we immediately have that switching the orientation of a surface multiplies the Poisson structure (the bi-vector) by $-1$.
\end{commnt}

We conclude this section by stating explicitly how the Poisson structure distinguishes the moduli of surface group representations where the algebraic structure alone cannot.

\begin{theorem}Let $\Sigma$ and $\Sigma'$ be compact, connected, orientable surfaces with boundary, and let $\X(\Sigma)$ and $\X(\Sigma')$ be the moduli of surface group representations into $\mathrm{SL}(m,\C)$.  Then:
\begin{enumerate}
\item $\X(\Sigma)$ and $\X(\Sigma')$ are isomorphic as varieties if and only if $\chi (\Sigma)=\chi (\Sigma')$
\item $\X(\Sigma)$ and $\X(\Sigma')$ are isomorphic as Poisson varieties if and only if $\Sigma$ and $\Sigma'$ are homeomorphic.
\end{enumerate}
\end{theorem}

\begin{proof}
$\chi (\Sigma)=\chi (\Sigma')$ if and only if the ranks of their fundamental groups are equal which implies $\X(\Sigma)$ and $\X(\Sigma')$ are isomorphic as varieties.  On the other hand, if they are 
isomorphic as varieties then they must share the same transcendence degree.  The transcendence degree is computed to be $-\chi(\Sigma)(m^2-1)$ in general.  Consequently, they have the same Euler 
characteristic. 

Observe $\Sigma$ and $\Sigma'$ are homeomorphic if and only if they have the same number of boundaries and the same genus.  Their moduli are equivalent as Poisson varieties if and only if they are isomorphic as varieties and there is a Lie algebra isomorphism of their coordinate rings.  With the same boundaries and genus, the moduli are certainly isomorphic by our previous work.  Moreover, there is certainly not a Lie algebra isomorphism if the Casimir subalgebras are different.  Since the bracket sums over intersections of cycles and in a surface the only cycles which necessarily are parallel to all other homotopy classes of cycles are the boundaries, the number of boundaries $n$ gives $n(m-1)$ generically independent generators of the Casimir subalgebra (since the boundary mapping is generically submersive).  Consequently, if $\X(\Sigma)$ and $\X(\Sigma')$ are isomorphic as Poisson varieties and $\Sigma$ and $\Sigma'$ are not homeomorphic we get a contradiction.  Equivalently, we have that the symplectic leaves have generic dimensions in terms of the number of boundaries.  So again, if the surfaces are not homeomorphic we have a contradiction since the Poisson structures are equivalent if the corresponding symplectic foliations are equivalent.  On the other hand, if $\Sigma$ and $\Sigma'$ are homeomorphic the intersection pairing on cohomology will be equivalent.  Since the Poisson bracket is in terms of intersections the moduli will have equivalent brackets.  This will allow for a Lie algebra isomorphism of their coordinate rings.
\end{proof}

\section{Structure of $\C[\G^{\times 2}\aq \G]$}\label{lawtonresults}
We now restrict our attention to $m=3$; that is, to $\SL$ character varieties.  The reason for this is that an explicit description of the Poisson bracket requires knowledge of the variety.  
However, an explicit description of the algebraic structure is not known in general.  In fact, the only general formulation is for $\mathrm{SL}(2,\C)$.  In an upcoming paper we explore more 
generally the Poisson structure for $\mathrm{SL}(2,\C)$ character varieties.  For $\SL$ character varieties the ideal of relations defining the character variety is only known for rank $1$ and $2$ 
free groups.  In this section we review results in \cite{L} that we will need to characterize the Poisson structures on $\SL$ character varieties of Euler characteristic $-1$ surfaces.

\subsection{Minimal Generators} From Section $\ref{repchar}$ we know that $\C[\X]$ is generated by traces of words in generic matrices of word length no greater than $6$.  The following lemma may be 
found in \cite{L}.

\begin{lemma}\label{gens}
$\C [\G\times \G]^\G$ is minimally generated by
\begin{align*}
&\tr{\xb_1},\ \ \tr{\xb_2},\ \ \tr{\xb_1\xb_2},\ \ \tr{\xb_1\xb_2^{-1}},\ \ \tr{\xb_1^{-1}},\\
&\tr{\xb_2^{-1}},\ \ \tr{\xb_1^{-1}\xb_2^{-1}},\ \ \tr{\xb_1^{-1}\xb_2},\ \ \tr{\xb_1\xb_2\xb_1^{-1}\xb_2^{-1}}.
\end{align*}
\end{lemma}

Since we are working with invariants of pairs of $3\times 3$ unimodular matrices the Krull dimension is $8$.  With $9$ generators, we expect one relation.

\subsection{Hyper-Surface in $\C^9$}  With a minimal generating set, it remains to describe the relations in terms of these generators.

Let $$\overline{R}=\C[t_{(1)},t_{(-1)},t_{(2)},t_{(-2)},t_{(3)},t_{(-3)},t_{(4)},t_{(-4)},
t_{(5)},t_{(-5)}]$$ be the complex polynomial ring freely generated by $\{t_{(\pm i)},\ 1\leq i\leq 5\},$ and let
$$R=\C[t_{(1)},t_{(-1)},t_{(2)},t_{(-2)},t_{(3)},t_{(-3)},t_{(4)},t_{(-4)}]$$ be its subring generated
by $\{t_{(\pm i)},\  1\leq i\leq 4\},$  so
$\overline{R}=R[t_{(5)},t_{(-5)}].$
Define the following ring homomorphism,
$$R[t_{(5)},t_{(-5)}]\stackrel{ \Pi}{\longrightarrow} \C[\G\times \G]^\G$$ by

\begin{center}
\begin{tabular}{ll}
$t_{(1)}\mapsto\tr{\xb_1}$ & $t_{(-1)}\mapsto\tr{\xb_1^{-1}}$\\
$t_{(2)}\mapsto\tr{\xb_2}$& $t_{(-2)}\mapsto\tr{\xb_2^{-1}}$\\
$t_{(3)}\mapsto\tr{\xb_1\xb_2}$& $t_{(-3)}\mapsto\tr{\xb_1^{-1}\xb_2^{-1}}$\\
$t_{(4)}\mapsto\tr{\xb_1\xb_2^{-1}}$& $t_{(-4)}\mapsto\tr{\xb_1^{-1}\xb_2}$\\
$t_{(5)}\mapsto\tr{\xb_1\xb_2\xb_1^{-1}\xb_2^{-1}}$& $t_{(-5)}\mapsto\tr{\xb_2\xb_1\xb_2^{-1}\xb_1^{-1}}$.
\end{tabular}
\end{center}

It follows from Lemma \ref{gens} that $$\C[\X]\cong R[t_{(5)}, t_{(-5)}]/\ker(\Pi).$$ In other words, $\Pi$ is a
surjective algebra morphism.
We define $$P=t_{(1)}t_{(-1)}t_{(2)}t_{(-2)}-t_{(1)}t_{(2)}t_{(-3)}-t_{(-1)}t_{(-2)}t_{(3)}
-t_{(1)}t_{(-2)}t_{(-4)}-t_{(-1)}t_{(2)}t_{(4)}$$
$$+t_{(1)}t_{(-1)}+t_{(2)}t_{(-2)}+t_{(3)}t_{(-3)}+t_{(4)}t_{(-4)}-3,$$ and so $P\in R$.
It is shown in \cite{L} that $$P-(t_{(5)}+t_{(-5)}) \in \ker(\Pi).$$

Hence it follows that the composite map $$R[t_{(5)}]\hookrightarrow R[t_{(5)}, t_{(-5)}]\twoheadrightarrow R[t_{(5)},
t_{(-5)}]/\ker(\Pi),$$ is an epimorphism.  Let $I$ be the kernel of this composite map. In \cite{L} it is then shown that there exists $Q\in
R$ so $Q-t_{(5)}t_{(-5)}\in \ker(\Pi)$ as well.

As a consequence we have the following lemma describing the ideal of relations defining the character variety.

\begin{lemma}\label{ideallemma}
$I$ is principally generated by the polynomial
\begin{equation}\label{idealequ}
t_{(5)}^2-Pt_{(5)}+Q.
\end{equation}
\end{lemma}

The form of $Q$ is shown to be
\begin{align}\label{q}
Q=&9-6t_{(1)}t_{(-1)} -6t_{(2)}t_{(-2)} -6t_{(3)}t_{(-3)} -6t_{(4)}t_{(-4)}+t_{(1)}^3 +t_{(2)}^3 +t_{(3)}^3 +t_{(4)}^3\nonumber\\
&+t_{(-1)}^3 +t_{(-2)}^3 +t_{(-3)}^3 +t_{(-4)}^3 -3t_{(-4)}t_{(-3)}t_{(-1)} -3t_{(4)}t_{(3)}t_{(1)} -\nonumber\\
&3t_{(-4)}t_{(2)}t_{(3)} -3t_{(4)}t_{(-2)}t_{(-3)}+3t_{(-4)}t_{(-2)}t_{(1)} +3t_{(4)}t_{(2)}t_{(-1)}+ \nonumber\\
&3t_{(1)}t_{(2)}t_{(-3)} +3t_{(-1)}t_{(-2)}t_{(3)}+t_{(-2)}t_{(-1)}t_{(2)}t_{(1)}+t_{(-3)}t_{(-2)}t_{(3)}t_{(2)} +\nonumber\\
&t_{(-4)}t_{(-1)}t_{(4)}t_{(1)} +t_{(-4)}t_{(-2)}t_{(4)}t_{(2)} +t_{(-3)}t_{(-1)}t_{(3)}t_{(1)}+\nonumber \\
&t_{(-3)}t_{(-4)}t_{(3)}t_{(4)}+t_{(-4)}^2t_{(-3)}t_{(-2)}  +t_{(4)}^2t_{(3)}t_{(2)} +t_{(-1)}^2t_{(-2)}t_{(-4)} +t_{(1)}^2t_{(2)}t_{(4)}+\nonumber\\
&t_{(1)}t_{(-2)}^2t_{(-3)} +t_{(-1)}t_{(2)}^2t_{(3)} +t_{(-4)}t_{(-3)}t_{(1)}^2 +t_{(4)}t_{(3)}t_{(-1)}^2 +\nonumber\\
&t_{(-4)}t_{(2)}t_{(-3)}^2 +t_{(4)}t_{(-2)}t_{(3)}^2 +t_{(-1)}^2t_{(-3)}t_{(2)} +t_{(1)}^2t_{(3)}t_{(-2)} +\nonumber\\
&t_{(-4)}t_{(1)}t_{(2)}^2 +t_{(4)}t_{(-1)}t_{(-2)}^2+t_{(-4)}t_{(3)}t_{(-2)}^2 +t_{(4)}t_{(-3)}t_{(2)}^2 +\nonumber\\
&t_{(1)}t_{(3)}t_{(-4)}^2 +t_{(-1)}t_{(-3)}t_{(4)}^2 +t_{(-1)}t_{(-4)}t_{(3)}^2+t_{(1)}t_{(4)}t_{(-3)}^2-2t_{(-3)}^2t_{(-2)}t_{(-1)}-\nonumber\\
&2t_{(3)}^2t_{(2)}t_{(1)} -2t_{(-4)}^2t_{(-1)}t_{(2)} -2t_{(4)}^2t_{(1)}t_{(-2)}+t_{(-1)}^2t_{(-2)}^2t_{(-3)}+t_{(1)}^2t_{(2)}^2t_{(3)}+\nonumber\\
&t_{(-4)}t_{(-1)}^2t_{(2)}^2+t_{(4)}t_{(1)}^2t_{(-2)}^2-t_{(-4)}t_{(-2)}^2t_{(2)}t_{(1)} -t_{(4)}t_{(2)}^2t_{(-2)}t_{(-1)}-\nonumber\\
&t_{(-3)}t_{(1)}^2t_{(-1)}t_{(2)}-t_{(3)}t_{(-1)}^2t_{(1)}t_{(-2)}- t_{(-3)}t_{(2)}^2t_{(-2)}t_{(1)} -t_{(3)}t_{(-2)}^2t_{(2)}t_{(-1)}-\nonumber\\
&t_{(-4)}t_{(-2)}t_{(-1)}t_{(1)}^2 -t_{(4)} t_{(2)}t_{(1)}t_{(-1)}^2-t_{(-1)}t_{(-2)}^3t_{(1)}-t_{(-1)}t_{(2)}^3 t_{(1)} -\nonumber\\
&t_{(-1)}^3t_{(-2)}t_{(2)}-t_{(1)}^3t_{(-2)}t_{(2)}-t_{(-4)}t_{(-3)}t_{(-2)}t_{(-1)}t_{(2)}-t_{(4)}t_{(3)}t_{(2)}t_{(1)}t_{(-2)}-\nonumber\\
&t_{(-1)}t_{(1)}t_{(2)}t_{(-4)}t_{(3)} -t_{(-1)}t_{(1)}t_{(-2)}t_{(4)}t_{(-3)}+ t_{(-2)}t_{(-1)}^2t_{(1)}^2t_{(2)} +t_{(-1)}t_{(-2)}^2t_{(2)}^2t_{(1)}.
\end{align}

Lemmas \ref{gens} and \ref{ideallemma} together imply 
\begin{theorem}\label{ranktwo}
$\G^{\times 2}\aq\G$ is isomorphic to an affine degree $6$ hyper-surface in $\C^9$, which maps onto $\C^8$ generically $2$-to-$1$.
\end{theorem}

With the generators and relations in hand, we now briefly describe symmetry within the character variety coming from the outer automorphisms of the surface group.  

\subsection{Outer Automorphisms}
Given any $\alpha\in \mathrm{Aut}(\F_2)$, we define $a_\alpha \in \mathrm{End}(\C[\X])$ by extending the following mapping
$$a_\alpha(\tr{\wb})=\tr{\alpha(\wb)}.$$  If $\alpha \in \mathrm{Inn}(\F_2)$, then there exists $\mathtt{u}\in \F_2$ so for all
$\mathtt{w}\in \F_2$, $$\alpha(\mathtt{w})=\mathtt{uwu}^{-1},$$ which implies
$$a_\alpha(\tr{\wb)}=\tr{\ub\wb\ub^{-1}}=\tr{\wb}.$$  

Thus $\ot(\F_2)$ acts on $\C[\X]$.  By results in \cite{MKS}, $\ot(\F_2)$ is generated by the following mappings
\begin{align}
\mathfrak{t}&=\left\{
\begin{array}{l}
\xt_1\mapsto \xt_2\\
\xt_2\mapsto \xt_1
\end{array}\right.\\
\mathfrak{i}_1&=\left\{
\begin{array}{l}
\xt_1\mapsto \xt_1^{-1}\\
\xt_2\mapsto \xt_2
\end{array}\right.\\
\mathfrak{n}&=\left\{
\begin{array}{l}
\xt_1\mapsto \xt_1\xt_2\\
\xt_2\mapsto \xt_2
\end{array}\right.
\end{align}

More generally for any $\alpha \in \ot(\F_r)$, $\tr{\alpha(\wb)}$ is a polynomial in traces of words; consequently $\ot(\F_r)$ acts by polynomials.

Let $\mathfrak{D}$ be the subgroup generated by $\mathfrak{t}$ and $\mathfrak{i}_1$, and let $\C\mathfrak{D}$ be the corresponding group ring.
Then $\C[\X]$ is a $\C\mathfrak{D}$-module.  It is not hard to work out that $\mathfrak{D}$ is isomorphic to the order 8 dihedral group $D_4$.  Moreover, the projection from Theorem \ref{ranktwo} is 
$\mathfrak{D}$-equivariant.

Using this action, one can establish the following succinct expressions for the polynomial relations $P$ and $Q$ (see \cite{L}).

\begin{corollary}
In $\C\mathfrak{D}$ define $\mathbb{S}_\mathfrak{D}={\displaystyle\sum_{\sigma\in \mathfrak{D}}\sigma.}$  Then  
$P=\mathbb{S}_\mathfrak{D}(p)-3$ and $Q=\mathbb{S}_\mathfrak{D}(q)+9$ where $p$ and $q$ are given by:
\begin{align*}
p=\frac{1}{8}\big(t_{(1)}t_{(-1)}t_{(2)}t_{(-2)}-4t_{(1)}t_{(-2)}t_{(-4)}+2t_{(1)}t_{(-1)}+&2t_{(3)}t_{(-3)}\big)\\
q=\frac{1}{8}\big(2t_{(-2)}t_{(-1)}^2 t_{(1)}^2t_{(2)}+4t_{(1)}^2t_{(2)}^2t_{(3)}-4t_{(1)}^3t_{(-2)}t_{(2)}&-8t_{(-4)}t_{(-2)}t_{(-1)}t_{(1)}^2-\\
4t_{(4)}t_{(3)}t_{(2)}t_{(1)}t_{(-2)}+8t_{(1)}t_{(3)}t_{(-4)}^2+8t_{(-4)}t_{(1)}t_{(2)}^2& -8t_{(3)}^2 t_{(2)}t_{(1)}+\\ 
4t_{(4)}t_{(-3)}t_{(2)}^2+t_{(-2)}t_{(-1)}t_{(2)}t_{(1)}+t_{(-3)}t_{(-4)}t_{(3)}t_{(4)}+&4t_{(-3)}t_{(-1)}t_{(3)}t_{(1)} +\\
4t_{(1)}^3 +4t_{(3)}^3+12t_{(-4)} t_{(-2)}t_{(1)}-12t_{(-4)}t_{(2)}t_{(3)}-&12t_{(1)}t_{(-1)} -12t_{(3)}t_{(-3)}\big).
\end{align*}
\end{corollary}

\section{Poisson Structure of a Trinion}\label{threehole}
We now can derive the explicit forms of Poisson brackets for Euler characteristic $-1$ surfaces; that is, surfaces with fundamental groups free of rank $2$.  We begin with the three-holed sphere.

\begin{theorem}\label{threeholebracket}
Let $\X$ be the relative character variety of $S=\Sigma_{3,0}$.  Then the Poisson bracket has the following properties which determines all pairing of elements in $\C[\X]$:   $t_{(\pm 1)}, \ti{\pm 2}, 
\ti{\pm 3}$ are Casimirs,
\begin{align}\mathfrak{a}_{4,-4}&=\{t_{(4)},t_{(-4)}\}=P-2t_{(5)}, \text{ and}\label{t4formula}\\
\mathfrak{a}_{\pm 4,5}&=\{t_{(\pm 4)},t_{(5)}\}=\frac{t_{(5)}\{t_{(\pm 4)},P\}-\{t_{(\pm 4)},Q\}}{\{t_{(-4)},t_{(4)}\}}=\pm \frac{\partial(Q-\ti{5}P)}{\partial \ti{\mp 
4}}.\label{t5formula}\end{align}
\end{theorem}

For any Poisson bracket there exists an exterior bi-vector field whose restriction to symplectic leaves corresponds to the symplectic form.  Denote this bi-vector by $\mathfrak{a}$, and let $f,g\in \C[\X]$.  Then 
with respect to interior multiplication $\{f,g\}=\mathfrak{a}\cdot df\otimes dg.$  In local coordinates $(z_1,...,z_k)$ it takes the form 
$$\mathfrak{a}=\sum_{i,j}\mathfrak{a}_{i,j}\frac{\partial}{\partial z_i}\land 
\frac{\partial}{\partial z_j}$$ and so 
\begin{align*}
\{f,g\}&=\sum_{i,j}\left(\mathfrak{a}_{i,j}\frac{\partial}{\partial z_i}\land \frac{\partial}{\partial z_j}\right)\cdot\left(\frac{\partial f}{\partial z_i} dz_i \otimes \frac{\partial g}{\partial 
z_j}dz_j\right)\\
&=\sum_{i,j}\mathfrak{a}_{i,j}\left(\frac{\partial f}{\partial z_i}\frac{\partial g}{\partial z_j}-
\frac{\partial f}{\partial z_j}\frac{\partial g}{\partial z_i}\right).
\end{align*}

Theorem \ref{threeholebracket}, the calculations derived in its proof (see below), and observing the symmetry between $\{\ti{4},\ti{5}\}$ and $\{\ti{-4},\ti{5}\}$ allows for a succinct expression of 
the Poisson bi-vector field and proves the following corollary.

\begin{corollary}
$\mathfrak{a}_{-4,5}=-\mathfrak{i}(\mathfrak{a}_{4,5})$ where $\mathfrak{i}=\mathfrak{i}_1\mathfrak{t}\mathfrak{i}_1\mathfrak{t}$ is the outer automorphism $\xt_i\mapsto \xt_i^{-1}$, and the Poisson bi-vector field for the 
relative $\SL$-character variety of $\Sigma_{3,0}$ is 
$$(P-2\ti{5})\frac{\partial}{\partial \ti{4}}\land \frac{\partial}{\partial \ti{-4}}
+(1-\mathfrak{i})\left(\frac{\partial (Q-\ti{5}P)}{\partial \ti{-
4}}\frac{\partial}{\partial \ti{4}}\land \frac{\partial}{\partial \ti{5}}\right).$$
\end{corollary}

\begin{proof}[Proof of Theorem \ref{threeholebracket}]  Since a Poisson bracket is a bilinear, anti-commutative derivation, it is completely determined once it is formulated on the generators of 
$\C[\X]$.  

We present the fundamental group of $S=\Sigma_{3,0}$ as $$\pi_1(S)=\{\xt_1,\xt_2,\xt_3\ | \ \xt_3 \xt_2 \xt_1=1\},$$ 
so $\xt_3=\xt_1^{-1}\xt_2^{-1}$.  Hence $\pi_1(S)$ is free of rank $2$.

\begin{figure}[!ht]
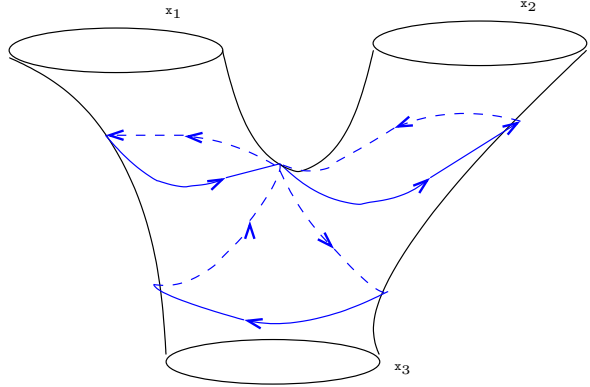

\begin{center}
\include{pantspi1}
\caption{Presentation of $\pi_1(\Sigma_{3,0},*)$}\label{pantspi1fig}
\end{center}
\end{figure}

The boundary curves in $S$ are the words $\xt_1$, $\xt_2$, and $(\xt_2\xt_1)^{-1}$, which are disjoint in the 
surface.  The sum in formula $\eqref{bracket}$ is taken over intersections, and is well-defined on homotopy 
classes.  So the trace of words corresponding to disjoint curves Poisson commute; that is, they 
are Casimirs.  Hence $t_{(\pm i)}$ are Casimirs, for $1\leq i\leq 3$, since they correspond to traces of boundary curves 
(and their inverses) in $S$.  

Using the derivation property and the identity $t_{(5)}^2-Pt_{(5)}+Q=0$, we deduce:
$$t_{(5)} \{t_{(\pm 4)},P\}+P \{t_{(\pm 4)},t_{(5)}\}-\{t_{(\pm 4)},Q\}=\{t_{(\pm 4)},t_{(5)}^2\}=2t_{(5)}\{t_{(\pm 
4)},t_{(5)}\}.$$  Hence $$(2t_{(5)}-P) \{t_{(\pm 4)},t_{(5)}\}=t_{(5)} \{t_{(\pm 4)},P\}- \{t_{(\pm 4)},Q\}.$$  So $\eqref{t5formula}$ follows from $\eqref{t4formula}$.  

Now assuming $\eqref{t4formula}$ and subsequently using the explicit 
expressions of $P$ and $Q$ given in Section $\ref{lawtonresults}$, we further derive explicit expressions for $\eqref{t5formula}$ as follows.

\begin{align*}
\{t_{(4)},P\}=&(P-2t_{(5)})(t_{(4)}-t_{(1)}t_{(-2)})\\
\{t_{(4)},Q\}=&(P-2t_{(5)})(-6t_{(4)}+3t_{(-4)}^2-3t_{(-1)}t_{(-3)}-3t_{(2)}t_{(3)}+3t_{(1)}t_{(-2)}+\\
&t_{(1)}t_{(-1)}t_{(4)}+t_{(2)}t_{(-2)}t_{(4)}+t_{(3)}t_{(-3)}t_{(4)}+t_{(-1)}^2t_{(-2)}+t_{(1)}^2t_{(-3)}+\\
&t_{(2)}t_{(-3)}^2+t_{(1)}t_{(2)}^2+t_{(3)}t_{(-2)}^2+t_{(-1)}t_{(3)}^2 
+t_{(-1)}^2t_{(2)}^2-t_{(1)}t_{(-1)}t_{(2)}t_{(3)}-\\
&t_{(-3)}t_{(-2)}t_{(-1)}t_{(2)}-t_{(1)}t_{(2)}t_{(-2)}^2-t_{(-2)}t_{(-1)}t_{(1)}^2+2t_{(1)}t_{(3)}t_{(-4)}+\\
&2t_{(-2)}t_{(-3)}t_{(-4)}-4t_{(-1)}t_{(2)}t_{(-4)})
\end{align*}
\begin{align*}
\{t_{(-4)},P\}=&(2t_{(5)}-P)(t_{(-4)}-t_{(-1)}t_{(2)})\\
\{t_{(-4)},Q\}=&(2t_{(5)}-P)(-6t_{(-4)}+3t_{(4)}^2-3t_{(1)}t_{(3)}-3t_{(-2)}t_{(-3)}+3t_{(-1)}t_{(2)}+\\
&t_{(1)}t_{(-1)}t_{(-4)}+t_{(2)}t_{(-2)}t_{(-4)}+t_{(3)}t_{(-3)}t_{(-4)}+t_{(1)}^2t_{(2)}+t_{(-1)}^2t_{(3)}+\\
&t_{(-2)}t_{(3)}^2+t_{(-1)}t_{(-2)}^2+t_{(-3)}t_{(2)}^2 
+t_{(1)}t_{(-3)}^2+t_{(1)}^2t_{(-2)}^2-t_{(1)}t_{(-1)}t_{(-2)}t_{(-3)}-\\
&t_{(3)}t_{(-2)}t_{(1)}t_{(2)}  -t_{(-1)}t_{(-2)}t_{(2)}^2 -t_{(2)}t_{(1)}t_{(-1)}^2+2t_{(-1)}t_{(-3)}t_{(4)}+\\
&2t_{(2)}t_{(3)}t_{(4)}-4t_{(1)}t_{(-2)}t_{(4)}) 
\end{align*}
and so 
\begin{align*}
\{t_{(4)},t_{(5)}\}=&t_{(4)}\big(t_{(1)}t_{(-1)}+t_{(2)}t_{(-2)}+t_{(3)}t_{(-3)}-t_{(5)}-6\big)+\\
&t_{(-4)}\big(2t_{(1)}t_{(3)}+2t_{(-2)}t_{(-3)}-4t_{(-1)}t_{(2)}\big)+\\
&t_{(5)}t_{(1)}t_{(-2)}+3t_{(-4)}^2-3t_{(-1)}t_{(-3)}-3t_{(2)}t_{(3)}+3t_{(1)}t_{(-2)}+\\
&t_{(-1)}^2t_{(-2)} +t_{(1)}^2t_{(-3)}+t_{(2)}t_{(-3)}^2 +t_{(1)}t_{(2)}^2 +t_{(3)}t_{(-2)}^2+\\
&t_{(-1)}t_{(3)}^2 +t_{(-1)}^2t_{(2)}^2 -t_{(1)}t_{(-1)}t_{(2)}t_{(3)}-t_{(-3)}t_{(-2)}t_{(-1)}t_{(2)} -\\
&t_{(1)}t_{(2)}t_{(-2)}^2 -t_{(-2)}t_{(-1)}t_{(1)}^2
\end{align*}
\begin{align*}
\{t_{(-4)},t_{(5)}\}=&t_{(-4)}\big(t_{(5)}-t_{(-1)}t_{(1)}-t_{(2)}t_{(-2)}-t_{(3)}t_{(-3)}+6\big)+\\
&t_{(4)}\big(4t_{(1)}t_{(-2)}-2t_{(-1)}t_{(-3)}-2t_{(2)}t_{(3)}\big)-\\
&t_{(5)}t_{(-1)}t_{(2)}-3t_{(4)}^2+3t_{(1)}t_{(3)}+3t_{(-2)}t_{(-3)}-3t_{(-1)}t_{(2)}-\\
&t_{(1)}^2t_{(2)} -t_{(-1)}^2t_{(3)}-t_{(-2)}t_{(3)}^2 -t_{(-1)}t_{(-2)}^2 -t_{(-3)}t_{(2)}^2-\\
&t_{(1)}t_{(-3)}^2 -t_{(1)}^2t_{(-2)}^2 +t_{(-1)}t_{(1)}t_{(-2)}t_{(-3)}+t_{(3)}t_{(2)}t_{(1)}t_{(-2)} +\\
&t_{(-1)}t_{(-2)}t_{(2)}^2 +t_{(2)}t_{(1)}t_{(-1)}^2.
\end{align*}

It remains to compute $ \{t_{(4)},t_{(-4)}\}$.  Following the results in \cite{G5}, we consider immersed closed curves freely homotopic to $\alpha=\xt_1\xt_2^{-1}$ and $\beta=\xt_2\xt_1^{-1}$ 
intersecting only at transverse double points.
  
\begin{figure}[!ht]
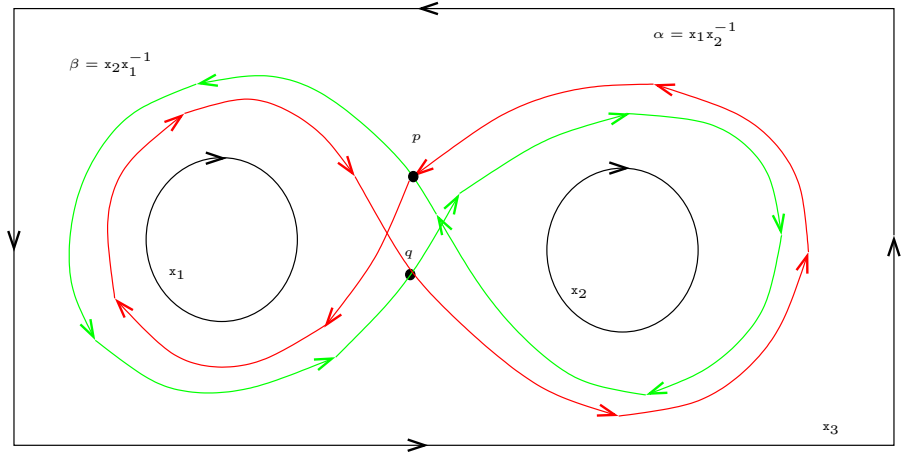

\begin{center}
\include{t4.new}
\caption{$\alpha$ and $\beta$ in $S$}\label{t4fig}
\end{center}
\end{figure}

Since $S$ is homotopic to a closed rectangle with two open disks removed, we depict all curves as in Figure $\ref{t4fig}$. 

We further let $\alpha_p$ and $\beta_p$ be the curves corresponding to $\alpha$ and $\beta$ based at the point $p$ in 
$\pi_1(S,p)$.

\begin{figure}[!ht]
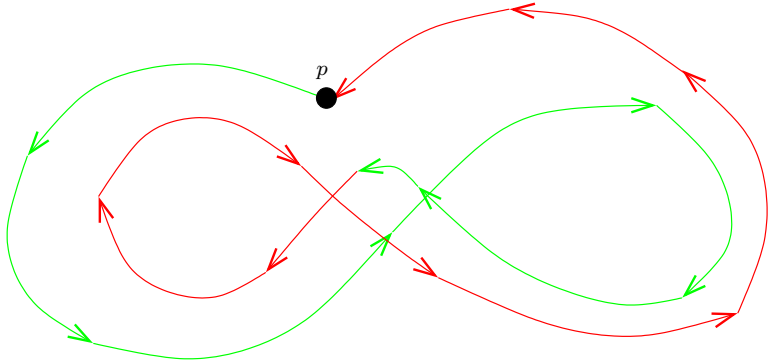

\begin{center}
\include{t4p.new}
\caption{$\alpha_p\beta_p=\xt_2^{-1}\xt_1\xt_2\xt_1^{-1}$}\label{t4pfig}
\end{center}
\end{figure}

Respectively, let $\alpha_q$ and $\beta_q$ be the corresponding curves in $\pi_1(S,q)$.

\begin{figure}[!ht]
\begin{center}
\include{t4q.new}
\caption{$\alpha_q\beta_q=\xt_1\xt_2^{-1}\xt_1^{-1}\xt_2$}\label{t4qfig}
\end{center}
\end{figure}

Calculating the oriented intersection number at $p$ and $q$ we find $\epsilon(p,\alpha,\beta)=-1$ and  
$\epsilon(q,\alpha,\beta)=1.$

\begin{figure}[!ht]
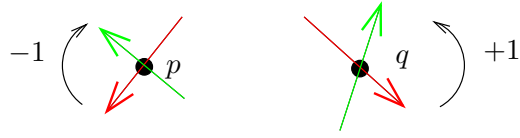

\begin{center}
\include{pqsign}
\caption{Intersection numbers at $p$ and $q$}\label{pqsignfig}
\end{center}
\end{figure}

Hence formula $\eqref{bracket}$ and Figures $\ref{t4pfig}, \ref{t4qfig},$ and $\ref{pqsignfig}$ give
\begin{align*}
\{t_{(4)},t_{(-4)}\}=&\{\mathrm{tr}(\rho(\alpha)),\mathrm{tr}(\rho(\beta))\}\\
=&\epsilon(p,\alpha,\beta)\big(\mathrm{tr}(\rho(\alpha_p\beta_p))-(1/3)\mathrm{tr}
(\rho(\alpha))\mathrm{tr}(\rho(\beta))\big)+\\
&\epsilon(q,\alpha,\beta)\big(\mathrm{tr}(\rho(\alpha_q\beta_q))-(1/3)\mathrm{tr}
(\rho(\alpha))\mathrm{tr}(\rho(\beta))\big)\\
=&-\mathrm{tr}(\rho(\alpha_p\beta_p))+\mathrm{tr}(\rho(\alpha_q\beta_q))\\
=&-\mathrm{tr}(\xb_2^{-1}\xb_1\xb_2\xb_1^{-1})+\mathrm{tr}(\xb_1\xb_2^{-1}\xb_1^{-1}\xb_2)\\
=&-t_{(5)}+t_{(-5)}\\
=&-t_{(5)}+(P-t_{(5)})=P-2t_{(5)}.
\end{align*}\qedhere\end{proof}

\begin{commnt}
Formula $\eqref{t5formula}$ can be derived in the same manner as we derived formula 
$\eqref{t4formula}$.  Doing so leads to the expression:
$$\mathrm{tr}(\xb_1\xb_2^{-1}\xb_1^{-1}\xb_2^{-1}\xb_1\xb_2)-
\mathrm{tr}(\xb_1\xb_2^{-1})+\mathrm{tr}(\xb_2^{-2}\xb_1^{2}\xb_2\xb_1^{-1})-
\mathrm{tr}(\xb_2^{-1}\xb_1\xb_2^{-1}\xb_1\xb_2\xb_1^{-1}).$$
Subsequently using polynomial trace relations to reduce these trace expressions to polynomials in 
$t_{(i)}$ for $1\leq |i|\leq 5$ has provided us with further verification of $\eqref{t5formula}$. 
\end{commnt}

\section{Poisson Structure on a One-Holed Torus}\label{onehole}

The only other orientable surface with Euler characteristic $-1$ is the one-holed torus; that is, $S^1\times S^1$ with a disk removed.  Let $T=\Sigma_{1,1}$ be the one-holed torus.  Its fundamental group is presentable as $$\pi_1(T)=\{\xt_1, \xt_2, \xt_3\ |\ \xt_1\xt_2\xt_1^{-1}\xt_2^{-1}\xt_3=1\}.$$

\begin{figure}[!ht]
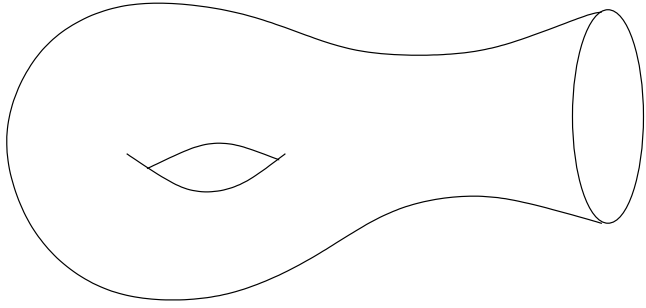

\begin{center}
\include{torus_uncut}
\caption{One-Holed Torus}\label{torus_uncut}
\end{center}
\end{figure}

With respect to this presentation the boundary  $\xt_3$ corresponds to the inverse of the word $\xt_1\xt_2\xt_1^{-1}\xt_2^{-1}$.  Consequently, the only
Casimir that is also a generator is $\ti{5}$ since it corresponds to $\tr{\xb_1\xb_2\xb_1^{-1}\xb_2^{-1}}$.  Thus the 81 pairings coming from the 9 generators is reduced to 64.  Anticommutativity 
reduces this number to 28.

Observe that $\{\ti{i},\ti{-i}\}=0$ since the corresponding curves are homotopic in $T$ to parallel curves and so have no intersection. Note that this is not always true in general; in the case of the 
three-holed sphere $\{\ti{4},\ti{-4}\}\not=0$.   We are left with 24 computations.

We first address the pairings which come from cycles of word length one and so only have one intersection.

Let $\alpha$ be homotopic to $\xt_1$ and $\beta$ be homotopic to $\xt_2$, as in Figure $\ref{torus:1}$

\begin{figure}[!ht]
\begin{center}
\include{torus}
\caption{$\alpha=\xt_1$ and $\beta=\xt_2$ in $T$}\label{torus:1}
\end{center}
\end{figure}

Let $\epsilon(p,\alpha, \beta)=\epsilon_p$ to simplify the notation.

Then 
\begin{align*}
\{\ti{1},\ti{2}\}&=\epsilon_p\big(\tr{\rho(\alpha_p\beta_p)}-(1/3)\tr{\rho(\alpha)}\tr{\rho(\beta)}\big)\\
                 &=\epsilon_p\big(\tr{\xb_1\xb_2}-(1/3)\tr{\xb_1}\tr{\xb_2}\big)\\
                 &=\tr{\xb_1\xb_2}-(1/3)\tr{\xb_1}\tr{\xb_2}\\
                 & =\ti{3}-(1/3)\ti{1}\ti{2}.
\end{align*}

Likewise, we can compute
\begin{align*}
\{\ti{-1}, \ti{2}\}&=-\ti{-4}+(1/3)\ti{-1}\ti{2},\\ 
\{\ti{1}, \ti{-2}\}&=-\ti{4}+(1/3)\ti{1}\ti{-2},\\
\{\ti{-1}, \ti{-2}\}&=\ti{-3}-(1/3)\ti{-1}\ti{-2}.  
\end{align*}

We can already see symmetry coming from $\ot(\F_2)$.  In $\mathfrak{D}\cong D_4$ define $\mathfrak{i}_2=\mathfrak{t}\mathfrak{i}_1\mathfrak{t}$; the mapping which sends $\xt_2\mapsto \xt_2^{-1}$.  Then

\begin{align*}
\{\ti{-1}, \ti{2}\}&=-\mathfrak{i}_1\{\ti{1},\ti{2}\}=\{\mathfrak{i}_1\ti{1},\mathfrak{i}_1\ti{2}\},\\
\{\ti{1}, \ti{-2}\}&=-\mathfrak{i}_2\{\ti{1},\ti{2}\}=\{\mathfrak{i}_2\ti{1},\mathfrak{i}_2\ti{2}\},\\
\{\ti{-1}, \ti{-2}\}&=\mathfrak{i}\{\ti{1},\ti{2}\}=\{\mathfrak{i}\ti{1},\mathfrak{i}\ti{2}\}.
\end{align*}

We are left with 20 computations.

Next we consider curves which again intersect only once, but at least one of them has word length two.  There will be two cases: either you have a repeated letter or you get a cancellation 
after cyclic reduction.  We demonstrate both cases here and then state the other pairings of these types.  First we consider the case of the paired cycles having a common letter and no letter between 
them cancels.

For instance, let $\alpha=\xt_1$ and $\beta=\xt_1\xt_2$ (both have the letter $\xt_1$ and when concatenated the word is cyclically reduced).  

\begin{figure}[!ht]
\begin{center}
\include{torus2}
\caption{$\alpha=\xt_1$ and $\beta=\xt_1\xt_2$ in $T$}\label{torus:2}
\end{center}
\end{figure}

Then
\begin{align*}
\{\ti{1},\ti{3}\}&=\epsilon_p\big(\tr{\rho(\alpha_p\beta_p)}-(1/3)\tr{\rho(\alpha)}\tr{\rho(\beta)}\big)\\
                 &=\tr{\xb_1^2\xb_2}-(1/3)\tr{\xb_1}\tr{\xb_1\xb_2}\\
                 & =-\ti{-1}\ti{2}+\ti{-4}+(2/3)\ti{1}\ti{3}.
\end{align*}

There are seven more cycles corresponding to minimal generators having a common letter and one intersection.  They all are computed like the one above.  Again there is apparent symmetry coming from 
$\mathfrak{D}$.
\begin{align*}
\{\ti{-1}, \ti{-3}\}&=-\ti{1}\ti{-2}+\ti{4}+(2/3)\ti{-1}\ti{-3}=\mathfrak{i}\{\ti{1},\ti{3}\}=\{\mathfrak{i}\ti{1},\mathfrak{i}\ti{3}\},\\
\{\ti{1}, \ti{4}\}&=\ti{-1}\ti{-2}-\ti{-3}-(2/3)\ti{1}\ti{4}=-\mathfrak{i}_2\{\ti{1},\ti{3}\}=\{\mathfrak{i}_2\ti{1},\mathfrak{i}_2\ti{3}\},\\
\{\ti{-1}, \ti{-4}\}&=\ti{1}\ti{2}-\ti{3}-(2/3)\ti{-1}\ti{-4}=-\mathfrak{i}_1\{\ti{1},\ti{3}\}=\{\mathfrak{i}_1\ti{1},\mathfrak{i}_1\ti{3}\},\\
\{\ti{2}, \ti{3}\}&=\ti{-2}\ti{1}-\ti{4}-(2/3)\ti{2}\ti{3}=-\mathfrak{t}\{\ti{1},\ti{3}\}=\{\mathfrak{t}\ti{1},\mathfrak{t}\ti{3}\},\\
\{\ti{-2}, \ti{-3}\}&=\ti{2}\ti{-1}-\ti{-4}-(2/3)\ti{-2}\ti{-3}=-\mathfrak{i}\mathfrak{t}\{\ti{1},\ti{3}\}=\{\mathfrak{i}\mathfrak{t}\ti{1},\mathfrak{i}\mathfrak{t}\ti{3}\},\\
\{\ti{2}, \ti{-4}\}&=-\ti{-2}\ti{-1}+\ti{-3}+(2/3)\ti{2}\ti{-4}=\mathfrak{i}_1\mathfrak{t}\{\ti{1},\ti{3}\}=\{\mathfrak{i}_1\mathfrak{t}\ti{1},\mathfrak{i}_1\mathfrak{t}\ti{3}\},\\
\{\ti{-2}, \ti{4}\}&=-\ti{2}\ti{1}+\ti{3}+(2/3)\ti{-2}\ti{4}=\mathfrak{i}_2\mathfrak{t}\{\ti{1},\ti{3}\}=\{\mathfrak{i}_2\mathfrak{t}\ti{1},\mathfrak{i}_2\mathfrak{t}\ti{3}\}.
\end{align*}

We are now left with 12 computations.

Now we consider the other case of only one intersection in the pair of cycles, where after concatenating there is some non-trivial (reduced word length) cyclic reduction in the word.  For instance, 
let $\alpha=\xt_1$ and $\beta=\xt_1^{-1}\xt_2^{-1}$.

\begin{figure}[!ht]
\begin{center}
\include{torus3}
\caption{$\alpha=\xt_1$ and $\beta=\xt_1^{-1}\xt_2^{-1}$ in $T$}\label{torus:3}
\end{center}
\end{figure}

Then
\begin{align*}
\{\ti{1},\ti{-3}\}&=\epsilon_p\big(\tr{\rho(\alpha_p\beta_p)}-(1/3)\tr{\rho(\alpha)}\tr{\rho(\beta)}\big)\\
                 &=(-1)\big(\tr{\xb_1\xb_2^{-1}\xb_1^{-1}}-(1/3)\tr{\xb_1}\tr{\xb_1^{-1}\xb_2^{-1}}\big)\\
                 & =-\ti{-2}+(1/3)\ti{1}\ti{-3}.
\end{align*}

There are seven more cycles corresponding to minimal generators having a canceled letter and one intersection.  They all are computed like the one above.  Again there is apparent symmetry coming from
$\mathfrak{D}$.

\begin{align*}
\{\ti{-1}, \ti{3}\}&=-\ti{2}+(1/3)\ti{-1}\ti{3}=\mathfrak{i}\{\ti{1},\ti{-3}\}=\{\mathfrak{i}\ti{1},\mathfrak{i}\ti{-3}\},\\
\{\ti{1}, \ti{-4}\}&=\ti{2}-(1/3)\ti{1}\ti{-4}=-\mathfrak{i}_2\{\ti{1},\ti{-3}\}=\{\mathfrak{i}_2\ti{1},\mathfrak{i}_2\ti{-3}\},\\
\{\ti{-1}, \ti{4}\}&=\ti{-2}-(1/3)\ti{-1}\ti{4}=-\mathfrak{i}_1\{\ti{1},\ti{-3}\}=\{\mathfrak{i}_1\ti{1},\mathfrak{i}_1\ti{-3}\},\\
\{\ti{2}, \ti{-3}\}&=\ti{-1}-(1/3)\ti{2}\ti{-3}=-\mathfrak{t}\{\ti{1},\ti{-3}\}=\{\mathfrak{t}\ti{1},\mathfrak{t}\ti{-3}\},\\
\{\ti{-2}, \ti{3}\}&=\ti{1}-(1/3)\ti{-2}\ti{3}=-\mathfrak{i}\mathfrak{t}\{\ti{1},\ti{-3}\}=\{\mathfrak{i}\mathfrak{t}\ti{1},\mathfrak{i}\mathfrak{t}\ti{-3}\},\\
\{\ti{2}, \ti{4}\}&=-\ti{1}+(1/3)\ti{2}\ti{4}=\mathfrak{i}_1\mathfrak{t}\{\ti{1},\ti{-3}\}=\{\mathfrak{i}_1\mathfrak{t}\ti{1},\mathfrak{i}_1\mathfrak{t}\ti{-3}\},\\
\{\ti{-2}, \ti{-4}\}&=-\ti{-1}+(1/3)\ti{-2}\ti{-4}=\mathfrak{i}_2\mathfrak{t}\{\ti{1},\ti{-3}\}=\{\mathfrak{i}_2\mathfrak{t}\ti{1},\mathfrak{i}_2\mathfrak{t}\ti{-3}\}.
\end{align*}

We are now left with 4 computations.

The last case to consider is when there are two intersections.  For instance, let $\alpha=\xt_1\xt_2$ and $\beta=\xt_1\xt_2^{-1}$.

\begin{figure}[!ht]
\begin{center}
\include{torus4}
\caption{$\alpha=\xt_1\xt_2$ and $\beta=\xt_1\xt_2^{-1}$ in $T$}\label{torus:4}
\end{center}
\end{figure}

Then
\begin{align*}
\{\ti{3},\ti{4}\}&=\epsilon_p\big(\tr{\rho(\alpha_p\beta_p)}-(1/3)\tr{\rho(\alpha)}\tr{\rho(\beta)}\big)\\
                 &\ +\epsilon_q\big(\tr{\rho(\alpha_q\beta_q)}-(1/3)\tr{\rho(\alpha)}\tr{\rho(\beta)}\big)\\
                 &=(-1)\big(\tr{\xb_2\xb_1\xb_2^{-1}\xb_1}-(1/3)\tr{\xb_1\xb_2}\tr{\xb_1\xb_2^{-1}}\big)\\
                 &\ +(-1)\big(\tr{\xb_2\xb_1\xb_1\xb_2^{-1}}-(1/3)\tr{\xb_1\xb_2}\tr{\xb_1\xb_2^{-1}}\big)\\
                 &=-\tr{\xb_1^2}-\tr{\xb_1\xb_2\xb_1\xb_2^{-1}}+(2/3)\tr{\xb_1\xb_2}\tr{\xb_1\xb_2^{-1}}\\
                 & =-\ti{1}^2+\ti{-1}-\ti{-4}\ti{-2}-\ti{2}\ti{-3}+\ti{-1}\ti{2}\ti{-2}-\frac{1}{3}\ti{3}\ti{4}
\end{align*}

The last identity is a consequence of trace reduction formulas found in \cite{L}.  

In a like manner one can compute the follow other brackets of this type:

\begin{align*}
\{\ti{-3}, \ti{-4}\}&=-\ti{-1}^2+\ti{1}-\ti{4}\ti{2}-\ti{-2}\ti{3}+\ti{1}\ti{-2}\ti{2}-\frac{1}{3}\ti{-3}\ti{-4}=\mathfrak{i}\{\ti{3},\ti{4}\}=\{\mathfrak{i}\ti{3},\mathfrak{i}\ti{4}\},\\
\{\ti{3}, \ti{-4}\}&=\ti{2}^2-\ti{-2}+\ti{4}\ti{-1}+\ti{1}\ti{-3}-\ti{-2}\ti{1}\ti{-1}+\frac{1}{3}\ti{3}\ti{-4}=-\mathfrak{t}\{\ti{3},\ti{4}\}=\{\mathfrak{t}\ti{3},\mathfrak{t}\ti{4}\},\\
\{\ti{-3}, \ti{4}\}&=\ti{-2}^2-\ti{2}+\ti{-4}\ti{1}+\ti{-1}\ti{3}-\ti{2}\ti{-1}\ti{1}+\frac{1}{3}\ti{-3}\ti{4}=-\mathfrak{i}\mathfrak{t}\{\ti{3},\ti{4}\}=\{\mathfrak{i}\mathfrak{t}\ti{3},\mathfrak{i}\mathfrak{t}\ti{4}\}.
\end{align*}

This finishes the bracket computations for the one-holed torus.  Observe that the symmetry operators consistently respect Poisson bracket up to sign in the same way for all minimal generators.  
In other words, given any element $\mathfrak{d}\in \mathfrak{D}$, $$\mathfrak{d}\left(\{\ti{i},\ti{j}\}\right)=\pm\{\mathfrak{d}(\ti{i}),\mathfrak{d}(\ti{j})\},$$ for all $i,j$.  This 
compounded with the fact that the bracket is an anti-commutative derivation establishes:  

\begin{corollary}
$\mathfrak{D}$ consists of Poisson homomorphisms and Poisson anti-homomorphisms.
\end{corollary}

Define the following elements of the group ring of $\mathfrak{D}\cong D_4$ acting on $\C[\X]$:

\begin{itemize}
\item $\Sigma_{1}=1+\mathfrak{i}-\mathfrak{i}_1-\mathfrak{i}_2$
\item $\Sigma_{2}=1+\mathfrak{i}-\mathfrak{t}-\mathfrak{i}\mathfrak{t}$
\end{itemize}
and let $\mathfrak{a}_{i,j}=\{\ti{i},\ti{j}\}$. \\

Note $$\frac{1}{2}\Sigma_1\Sigma_2=1+\mathfrak{i}-\mathfrak{i}_1-\mathfrak{i}_2-\mathfrak{t}-\mathfrak{i}\mathfrak{t}+\mathfrak{i}_1\mathfrak{t}+\mathfrak{i}_2\mathfrak{t}.$$

Then putting our work together from this section proves:

\begin{theorem}
The Poisson bi-vector field for the $\SL$-relative character variety of the one-holed torus is
\begin{align*}&\Sigma_{1}\bigg(\mathfrak{a}_{1,2}\frac{\partial}{\partial \ti{1}}\land \frac{\partial}{\partial \ti{2}}\bigg)
+\Sigma_2\left(\mathfrak{a}_{3,4}\frac{\partial}{\partial \ti{3}}\land \frac{\partial}{\partial \ti{4}}\right)\\
&+\frac{1}{2}\Sigma_1\Sigma_{2}\bigg(\mathfrak{a}_{1,3}\frac{\partial}{\partial \ti{1}}\land \frac{\partial}{\partial \ti{3}}
+\mathfrak{a}_{1,-3}\frac{\partial}{\partial \ti{1}}\land \frac{\partial}{\partial \ti{-3}}\bigg),\end{align*}
where: 
\begin{itemize}
\item $\mathfrak{a}_{1,2}=\ti{3}-\frac{1}{3}\ti{1}\ti{2}$
\item $\mathfrak{a}_{1,3}=\frac{2}{3}\ti{1}\ti{3}-\ti{-1}\ti{2}+\ti{-4}$
\item $\mathfrak{a}_{1,-3}=-\ti{-2}+\frac{1}{3}\ti{1}\ti{-3}$
\item $\mathfrak{a}_{3,4}=-\ti{1}^2+\ti{-1}-\ti{-4}\ti{-2}-\ti{2}\ti{-3}+\ti{-1}\ti{2}\ti{-2}-\frac{1}{3}\ti{3}\ti{4}$.
\end{itemize}

\end{theorem}

Observe that $\mathfrak{a}_{1,3}=\mathfrak{n}_{(2)}(\mathfrak{a}_{1,2})$ and $\mathfrak{a}_{1,-3}=-\mathfrak{n}_{(-2)}(\mathfrak{a}_{1,2})$ where $\mathfrak{n}_{(2)}$ is the outer automorphism that 
sends $\xt_1\mapsto \xt_1$ and $\xt_2\mapsto \xt_1\xt_2$ and where $\mathfrak{n}_{(-2)}$ is the outer automorphism that
sends $\xt_1\mapsto \xt_1$ and $\xt_2\mapsto \xt_1^{-1}\xt_2^{-1}$.  So up to outer automorphisms the bracket is determined by $\mathfrak{a}_{1,2}$ and $\mathfrak{a}_{3,4}$ alone.

\begin{commnt}
Since the boundary of the one-holed torus is interior to its 2-cell and corresponds to a boundary on the exterior of the 2-cell of the three-holed sphere, the orientations are reversed with respect to each other in the above computations.  We describe how these two structures relate to each other in an upcoming paper, and this observation is important to understanding a sign difference between the two structures. 
\end{commnt}

\end{document}